\documentclass[12pt]{amsart}

\usepackage{enumitem}
\usepackage{verbatim}
\usepackage{xcolor}
\usepackage{hyperref}
\usepackage{amssymb,amsmath,amsthm,amsfonts}
\usepackage[margin=1.3in]{geometry}
\usepackage{mathtools}
\usepackage{graphicx}
\usepackage{mathrsfs}

\usepackage{tikz, environ}
\usetikzlibrary{arrows}
\usetikzlibrary{shapes,positioning,intersections,quotes}
\usetikzlibrary{decorations.pathreplacing}

\allowdisplaybreaks

\newcommand{\im}{\text{Im}}
\newcommand{\ph}{\varphi}
\newcommand{\eps}{\varepsilon}
\newcommand{\QQ}{\mathbb{Q}}
\newcommand{\NN}{\mathbb{N}}
\newcommand{\ZZ}{\mathbb{Z}}
\newcommand{\RR}{\mathbb{R}}
\newcommand{\BBB}{\mathcal{B}}
\newcommand{\LLL}{\mathcal{L}}
\newcommand{\per}{\mathrm{per}}

\newtheorem{theorem}{Theorem}[section]
\newtheorem{lemma}[theorem]{Lemma}
\newtheorem{proposition}[theorem]{Proposition}
\newtheorem{corollary}[theorem]{Corollary}
\newtheorem{example}[theorem]{Example}

\theoremstyle{remark}
\newtheorem{remark}[theorem]{Remark}

\theoremstyle{definition}
\newtheorem{definition}[theorem]{Definition}
\newtheorem{question}[theorem]{Question}

\numberwithin{equation}{section}

\title{Topologically mixing suspension flows over shift spaces}
\author{Jason Day}
\address{Dept.\ of Mathematics, University of Houston, Houston, TX 77204}
\email{jjday@uh.edu}

\date{\today}

\thanks{The author was partially supported by NSF grant DMS-2154378.}
\dedicatory{Dedicated to the memory of Todd Fisher}

\begin{document}

	
	\begin{abstract}
		We establish necessary and sufficient conditions for suspension flows over certain families of shift spaces to be topologically mixing.  We also show the similarities and differences between this case and the smooth measure theoretic setting on a manifold.  Additionally, we show that the set of roof functions defined on a shift space that produce suspension flows that are not topologically mixing is dense in the set of all continuous roof functions.
	\end{abstract}
	
	\maketitle
	
	\section{Introduction}
	Suspension flows (or special flows) are comprised of two components:
	a discrete dynamical system $f$ on a base space X paired with a roof function $r \colon X \to (0,\infty)$. There is a natural relation between a suspension flow and a Poincar\'e map given by a transversal to the flow and the return map, see for instance \cite[\S 1.2]{fish-hass}.
	
	The dynamical property we will focus on in this paper is topological mixing.  A map $f \colon X \to X$ is \emph{topologically mixing} if for any two open sets $U, V \subset X$, there exists an $N \in \NN$, such that for all $n\geq N$ we have that $f^n(U) \cap V \neq \emptyset$. Similarly, a continuous time dynamical system $\ph^t \colon X \to X$ is said to be \emph{topologically mixing} if for any open sets $U, V \subset X$ there exists a $T \in \RR$ such that for all $ t\geq T$, we have that $\ph^t(U) \cap V \neq \emptyset$.
	
	This paper will address two questions about topologically mixing suspension flows when the dynamics in the base is a shift space.  
	\begin{enumerate}
		\item What properties of a roof function guarantee that a suspension flow is topologically mixing?
		\item Is the set of continuous roof functions that induce topologically mixing suspension flows open and dense?
	\end{enumerate}

	\subsection{Characterizing topological mixing}
	These two questions are motivated by a number of important results.  Two such results for nontrivial (not a single periodic orbit) basic sets for an Axiom A flow appear in the work of Bowen and Plante in \cite{mixing-bowen} and \cite{Plante}, respectively.  These and other results are summarized in \cite[\S6]{fish-hass}.
	
	\begin{theorem}[\cite{fish-hass}]\label{thm:fish-hass-incom}
		If $\Lambda$ is a locally maximal hyperbolic set of a flow $\ph^t$ on a smooth connected manifold, then the following are equivalent.
		\begin{enumerate}
			\item $\ph^t\big|_\Lambda$ is topologically mixing.
			\item The periodic points of $\Lambda$ are dense in $\Lambda$ and their strong stable and strong unstable manifolds are dense in $\Lambda$.
			\item $\ph^t\big|_\Lambda$ is transitive, and each open set contains a collection of periodic points with setwise incommensurate periods (see Definition \ref{def:commensurate}).\label{item:incommensurate}
		\end{enumerate}
	\end{theorem}
	If the flow is a transitive Anosov flow, then combining \cite[Theorem 1.8]{Plante} with Theorem \ref{thm:fish-hass-incom} shows that the flow is either topologically mixing or it can be expressed as a suspension flow with a constant roof function modulo a time change.  Theorem \ref{thm:fish-hass-incom} is similar to the results of Bowen in \cite{mixing-bowen72,mixing-bowen} that show that if $\Lambda$ is a basic set, then not topologically mixing is equivalent to being a suspension flow with a constant roof function.
	
	We can classify roof functions in the continuous setting by adapting the cohomological condition from the measure theoretic setting.
	
	\begin{definition}
		Two roof functions $r,s \colon X \to \RR$ are \textit{cohomologous} if there is a continuous function $g \colon X \to X$ such that
		\[
		r(x) - s(x) = g(f(x)) - g(x).
		\]
		The function $g$ is called a \emph{transfer function}.
	\end{definition}
	
	Informally, one may think of the dichotomy from Bowen and Plante for topological mixing as 
	\[
	\text{not topologically mixing } \iff \text{ roof is cohomologous to a constant.}
	\]
	It is always true that a suspension flow with a roof function that is cohomologous to a constant is not topologically mixing, but the converse is more difficult to establish.  It is clear that a necessary condition is that the base map must be transitive. Many textbooks give an example of a topologically mixing homeomorphism of a compact metric space and a constant roof function to show that a suspension flow over a topologically mixing base need not be topologically mixing, and statements are often made about this being a rare phenomenon. However, to our knowledge, this has not been quantified in the literature. We obtain partial results to these questions and list some open questions at the end of \S\ref{sec:main-results}.

	\begin{remark}
		It would be interesting to see if this dichotomy holds more generally for suspension flows over connected metric spaces without leveraging hyperbolicity or a manifold structure; however, we provide examples in \S\ref{sec:no-mix-ex} and \S\ref{sec:no-mix-incommensurate} of suspension flows over shift spaces where the converse fails.  That is, there are roof functions that are not cohomologous to a constant but produce suspension flows that are not topologically mixing.
	\end{remark}

	Although these examples give us a negative answer, we will show that an analogous result to Theorem \ref{thm:fish-hass-incom} holds for certain families of subshifts (see Theorem \ref{thm:mix-dich} and Theorem \ref{thm:mix-dich-beta}).  To establish a dichotomy for topological mixing, we do not classify roof functions by a cohomologous condition; rather, we leverage a condition similar to item \ref{item:incommensurate} of Theorem \ref{thm:fish-hass-incom}.  In particular, we will use this condition to make a connection between period lengths under the flow and the presence or absence of topological mixing.
	
	\begin{definition}\label{def:commensurate}
		A subset $P \subset \RR$ is said to be \emph{setwise commensurate} if there exists a $\delta > 0$ such that $P \subset \delta \ZZ = \{\delta n : n \in \ZZ \}$.   We say they are \emph{setwise incommensurate} otherwise.  That is, the ratio of two elements of $P$ is irrational, or equivalently, $P$ generates a dense subgroup of $\RR$.
	\end{definition}
	For the majority of this paper, we will usually be referring to the collection of orbit lengths of the periodic points of the flow when we are discussing setwise commensurate or incommensurate sets.
	
	Relative to the measure theoretic setting, there is a scarcity of results regarding mixing for suspensions in the topological setting.  There are more results pertaining to mixing for suspension flows in the measure theoretic case for a number of reasons. One reason is that any flow on a measure space without fixed points is isomorphic to a suspension flow \cite{CornfeldFominSinai}.  Additionally, in \cite{FMT07}, the authors point out that due to the work of Sinai, Ruelle, and Bowen in the 1970s the notions of topologically mixing and (measure theoretical) mixing are equivalent for hyperbolic basic sets for smooth flows.  Moreover, mixing in the measure theoretical setting carries important probabilistic meaning that is not afforded in the topological sense like decay of correlations.
	
	There have also been efforts to identify which roof functions produce mixing or non-mixing suspension flows in the measure theoretic setting.  Suspension flows over an interval exchange transformation under a roof function of bounded variation is also not mixing \cite{Katok}.  Suspension flows over a rotation of the circle under a roof function of bounded variation is not mixing \cite{Kochergin}. The work in \cite{Ravotti} shows that for a certain class of parabolic flows a dichotomy similar to \cite{mixing-bowen72} holds.
	
	\subsection{Prevalence of topological mixing}
	
	The prevalence of mixing basic sets for an Axiom A flow has also been studied.  We highlight two key results.
	\begin{theorem}\label{thm:mixing-Bowen}\cite{mixing-bowen}
		The set of Axiom A flows in the $C^\ell$ topology such that a nontrivial basic set is not mixing is of the first category for $1\leq \ell \leq\infty$.  Furthermore, if $\BBB$ is the subset of the $C^\ell$ Anosov flows on a compact connected differentiable manifold that have a global cross-section which is an infranilmanifold, then $\BBB$ is open in the set of $C^\ell$ Anosov flows and the set of mixing flows of $\BBB$ is open and dense in $\BBB$.
	\end{theorem}
	
	\begin{theorem}\cite{FMT07}
		There exists a $C^2$ open and $C^\ell$ dense (for each $2 \leq \ell \leq \infty$) subset of Axiom A flows for which each nontrivial basic set is mixing.  Furthermore, there exists a $C^1$ open and $C^\ell$ dense (for $1 \leq \ell \leq \infty$) subset of Axiom A flows such that each nontrivial attracting basic set for the flow is mixing.
	\end{theorem}
	
	Theorem \ref{thm:mixing-Bowen} was extended beyond the Axiom A case by \cite{AAB04}.  They consider $C^1$ robustly transitive flows, which are a natural generalization of the hyperbolic setting. 
	
	The results in \cite{Ravotti} for parabolic flows also establish that the set of roof functions that induce mixing flows is dense in the space of continuous roof functions.  We show that an open and dense property like this or similar to \cite{mixing-bowen} or \cite{FMT07} fails for suspensions over shift spaces (see Theorem \ref{thm:no-mix-dense}).
	
	\subsection*{Acknowledgments} 
	Todd Fisher was my master's thesis advisor from 2018 - 2020.  He was a thoughtful and kind mentor, and I am deeply grateful for all he did for me. Todd Fisher played a major role in the origins of this project, and this paper is motivated by Example \ref{ex:not-mixing}, which he constructed.  Unfortunately, he passed away before we could fully investigate the phenomenon that appears in this example. Given these circumstances, I will take responsibility for any inaccuracies that may appear in this paper.
	
	I would also like to thank the anonymous referees for their careful reading and their helpful comments and insights that improved the quality of the paper and led to the inclusion of Theorems \ref{thm:SFT-cohomologous} and \ref{thm:SFT-G-delta}.
	
	\section{Background}\label{sec:background}
	\subsection{Symbolic dynamics}
	Let $A$ be a finite (or countable) set, which we call the \emph{alphabet}.  The elements of $A$ are called \emph{symbols}.  We consider the set of all bi-infinite sequences whose terms come from $A$ and denote it as $A^\ZZ$.  A point $x \in A^\ZZ$ can be written as $x = \dots x_{-2}x_{-1}.x_0x_1x_2\dots$ where the ``."  indicates the location of the 0th term in the bi-infinite sequence.  
	
	We give the set $A^\ZZ$ the metric $d(x,y) = 2^{-\min \{|n| : x_n \neq y_n\}}$.  If $A$ is finite, then the set $A^\ZZ$ is a compact, totally disconnected set.  The dynamics we consider is given by the \emph{(left) shift map} $\sigma \colon A^\ZZ\to A^\ZZ$, which is a homeomorphism defined by $(\sigma x)_n = x_{n+1}$. 
	\begin{definition}
		A \emph{(two-sided) shift space} is a closed, $\sigma$-invariant set $X\subset A^\ZZ$.
	\end{definition}
	
	A \emph{word} is a finite string of symbols $w \in A^n$ for some $n$.  If $w \in A^n$, then we write $|w| = n$ to denote the \emph{length} of $w$.  Given $x \in  A^\ZZ$ and $i,j \in \ZZ$ with $i < j$, let $x_{[i,j]}$ denote the word $x_i x_{i+1} \cdots x_{j-1}x_j$, and $x_{[i,j)}$ denote the word $x_i x_{i+1} \dots x_{j-1}$.
	Given a shift space $X$, the \emph{language} of $X$ is
	\[
	\LLL := \bigcup_{n=0}^\infty \LLL_n,
	\quad\text{where } \LLL_n \coloneqq \{w\in A^n : w \text{ appears in some } x \in X \}.
	\]
	We will occasionally use the notation $\LLL_{\geq n} = \{w \in \LLL : |w| \geq n\}$.
	
	One can define transitivity and topological mixing for shift spaces in terms of the language.
	
	\begin{definition}\label{def:mixing-words}
		A shift space is \emph{transitive} if, for any $u,v \in \LLL$, there exists a word $w \in \LLL$ such that $uwv \in \LLL$. A shift space is \emph{topologically mixing} if, for any $u,v \in \LLL$, there exists an $N \in \NN$ such that for any $n\geq N$, there is a word $w \in \LLL_n$ such that $uwv \in \LLL$.
	\end{definition}
	
	\begin{definition}
		A word $v$ is a \emph{synchronizing word} if for every pair $u,w\in\LLL$ such that $uv \in \LLL$ and $vw \in \LLL$ we also have that $uvw \in \LLL$.
	\end{definition}
	Synchronizing words can be used to concatenate blocks of words together.  Indeed, if $v$ is a synchronizing word and $vuv,vwv \in \LLL$, then $vuvwv \in \LLL$.  We can construct many periodic points by leveraging transitivity and a synchronizing word.
	
	Requiring that the shift possess a synchronizing word is a weaker condition than the shift space being a subshift of finite type.  Thus, Theorem \ref{thm:mix-dich} applies to a broader class of bases than subshifts of finite type.  Shifts with specification, $S$-gap shifts, and irreducible sofic shifts are all examples of shift spaces that possess synchronizing words \cite{bertrand,fischer}.
	
	If $w \in \LLL$, we denote the two-sided cylinder of $w$ by $[w]$.  In particular, if $|w|= 2n+1$, then we define
	\[
	[w] = \{x \in X : x_{[-n,n]} = w\}.
	\]
	Similarly, if $|w| = 2n$, then
	\[
	[w] = \{x \in X : x_{[-n,n)} = w\}
	\]
	
	We will also let $w^m = w \dots w$ where $w$ is repeated $m$ times.  We let $w^\infty$ denote the one-sided sequence $ww\dots$, and we will let $\overline{w}$ denote the bi-infinite sequence of repeated $w$'s $\dots www \dots$.  
	
	Note that we will occasionally abuse notation and write $w_i$ to denote a word rather than the $i$th symbol of $w$.  We will make it clear whenever we do this.
	
	\subsection{Suspension flows}
	
	\begin{definition}\label{SuspensionDefinition}\cite[p. 21]{BrinStuck}
		Given a map $f \colon X\to X$ and a function $r \colon X \to (0,\infty)$, consider the quotient space
		\[M_r=\{(x,t) \in X\times (0,\infty) : 0 \leq t \leq r(x)\} / {\sim}\]
		where ${\sim}$ is the equivalence relation $(x,r(x))\sim(f(x),0)$.  The \emph{suspension flow} of $f$ with \emph{roof function} $r(x)$ is the flow $\ph^t \colon M_r\to M_r$ defined by $\ph^T(x,s)=(f^n(x),s')$, where $n$ and $s'$ satisfy 
		\[\sum_{j=0}^{n-1}r(f^j(x))+s'=T+s, \hspace{10mm} 0\leq s'\leq r(f^n(x)).\]
		One may also think of a suspension flow as the quotient by ${\sim}$ of the vertical flow on $X \times \RR$.
	\end{definition}
	
	Theorems \ref{thm:mix-dich} and \ref{thm:mix-dich-beta} require that the roof function $r$ satisfies the Walters property from \cite{Wal78}.  In our setting, the Walters property is stronger than continuity but weaker than  H\"older continuity.  We formulate the definition here for symbolic spaces.
	
	\begin{definition}\label{def:Walters}
		A function $r \colon X \to \RR$ satisfies the \emph{Walters property} if for all $\eps > 0$, there exists a $k \geq 0$ such that for any $n \in \NN$ and $x,y \in X$ satisfying $x_{[-k,n+k]} = y_{[-k,n+k]}$, it follows that 
		\[
		\left| \sum_{j=0}^{n-1} r(\sigma^j x) - r(\sigma^j y) \right| < \eps.
		\]
	\end{definition}
	
	Although we do not use a cohomologous condition to classify topological mixing properties, we highlight an important connection between suspension flows with cohomologous roof functions. 
	
	\begin{proposition}\cite{fish-hass}\label{prop:conjugate}
		Suppose $X$ is a compact metric space, and $r_1$ and $r_2$ are cohomologous by a transfer function $g \colon X\to\RR$.  If $\ph^t_{1}$ is the suspension flow under $r_1$ and $\ph^t_2$ is the suspension flow under $r_2$, then $\ph^t_{1}$ and $\ph^t_{2}$ are conjugate.
	\end{proposition}
	
	As a consequence of Proposition \ref{prop:conjugate} which asserts that cohomologous roof functions induce conjugate flows, determining when two roof functions are cohomologous becomes an important question to answer.  The following proposition provides a quick way to verify when two roof functions are not cohomologous.
	
	\begin{proposition}\label{prop:not-cohomologous}
		Let $r$ and $s$ be roof functions, and suppose $p\in X$ is a periodic point with period $\per(p)$ under $f$.  If
		\[\sum_{j=0}^{\per(p)-1}r(f^j(p))-s(f^j(p))\neq0\]
		then $r$ and $s$ are not cohomologous.
	\end{proposition}

	The following theorems are well-known results that apply to transitive subshifts of finite type. The Closing Lemma dates back to \cite{dA67}, but we have written it here to suit our setting.
	
	\begin{theorem}[Closing Lemma]\label{lem:closing}
		If $X$ is a transitive subshift of finite type, then for every $\eps>0$ there exists $\gamma > 0$ such that if $x \in X$ and $n \geq 0$ are such that $d(\sigma^n x,x) < \gamma$, then there exists $y \in X$ such that $\sigma^n y = y$ and $d(\sigma^j y, \sigma^j x) < \epsilon$ for all $0 \leq j < n$.
	\end{theorem}
	
	The Liv\v sic theorem was proved for H\"older functions in \cite{aL72}, but it is also true for functions satisfying the Walters property \cite[Theorem 4]{tB01}.  We have written it here in its more general form but have replaced the ``coboundary" terminology with the ``cohomologous" language that is relevant in this context.
	
	\begin{theorem}[Liv\v sic Theorem]\label{thm:Livsic}
		Let $X$ be a compact metric space, $f\colon X\to X$ a continuous map satisfying the closing lemma and possessing a point whose orbit is dense, and $r,s \colon X\to \RR$ are continuous functions satisfying the Walters property. Then $r$ and $s$ are cohomologous if and only if for every periodic point $x = f^p(x) \in X$, we have 
		\[
		\sum_{j=0}^{n-1} r(f^j (x)) - s(f^j (x)) = 0.
		\]
	\end{theorem}
	
	\section{Main results}\label{sec:main-results}
	
	We now state our main theorems and include some open questions at the end of this section.
	
	\begin{theorem}\label{thm:mix-dich}
		Suppose $X$ is a transitive shift space (from a finite or countable alphabet) with a synchronizing word $v$ and the roof function $r\colon X \to (0,\infty)$ satisfies the Walters property. The flow induced by the roof function $r$ is not topologically mixing if and only if there exists a $\delta>0$ such that for every periodic point $p \in [v]$, we have
		\begin{equation}\label{eqn:no-mix-orbit-length}
			\sum_{j=0}^{\per(p)-1}r(\sigma^jp) \in \delta \ZZ,
		\end{equation}
		where $\per(p)$ is the period of $p$ under the shift map $\sigma$.
	\end{theorem}
	The proof of Theorem \ref{thm:mix-dich} can be found in \S\ref{sec:mix-dict-pf}.
	
	\begin{remark}
		The statement in Theorem \ref{thm:mix-dich} is written as a complement to Theorem \ref{thm:fish-hass-incom}.  Note that the definition of setwise commensurate does not mean that period lengths are all pairwise commensurate.  For the roof function to satisfy equation \eqref{eqn:no-mix-orbit-length} of Theorem \ref{thm:mix-dich}, the periods must all be collectively integer multiples of $\delta$.  We show that this distinction is necessary with an example in \S\ref{sec:mixing-example}.
	\end{remark}
	
	If $X$ is a subshift of finite type, then we can say more about the roof function $r$.
	
	\begin{definition}
		A function $r\colon X \to \RR$ is \emph{locally constant} if, for every point $x \in X$, there is a neighborhood $U$ of $x$ such that $f|_U$ is constant.
	\end{definition}
	
	\begin{theorem}\label{thm:SFT-cohomologous}
		Suppose $X$ is a transitive subshift of finite type, and $r \colon X \to (0,\infty)$ is a roof function satisfying the Walters property.  If the suspension flow associated to $X$ and $r$ is not topologically mixing, then there exists a $\delta > 0$ such that $r$ is cohomologous to a locally constant roof function $s \colon X \to \delta\NN$.
	\end{theorem}
	
	The space of Walters functions over a shift space can be given a norm that makes it a separable Banach space.  We describe the norm in \S\ref{sec:SFT}, and direct the reader to \cite{tB01} for further information about the norm and the space of Walters functions.  Under the topology induced by this norm, we obtain the following result.
	
	\begin{theorem}\label{thm:SFT-G-delta}
		If $X$ is a transitive subshift of finite type, then the set of roof functions satisfying the Walters property that yield a topologically mixing suspension flow over $X$ is a $G_{\delta}$ dense subset in the set of all Walters roof functions defined on $X$.
	\end{theorem}
	
	The proofs of Theorems \ref{thm:SFT-cohomologous} and \ref{thm:SFT-G-delta} can be found in \S\ref{sec:SFT}

	A characterization similar to Theorem \ref{thm:mix-dich} holds for suspension flows over $\beta$-shifts.
	\begin{theorem}\label{thm:mix-dich-beta}
		Suppose $X$ is a $\beta$-shift and $r\colon X \to (0,\infty)$ satisfies the Walters property.  The flow induced by the roof function $r$ is not topologically mixing if and only if  there exists a $\delta>0$ such that for every periodic point $p \in [0]$, we have
		\begin{equation}\label{eqn:no-mix-orbit-length-beta}
			\sum_{j=0}^{\per(p)-1}r(\sigma^jp) \in \delta \ZZ,
		\end{equation}
		where $\per(p)$ is the period of $p$ under the shift map $\sigma$.
	\end{theorem}
	The proof of Theorem \ref{thm:mix-dich-beta} can be found in \S\ref{sec:beta-shifts}.
	
	If the roof function is locally constant, then we can say more about the topological mixing properties of the suspension for arbitrary shift spaces on a finite alphabet.
	
	\begin{theorem}\label{thm:loc-const-com-nonmixing}
		Suppose $r$ is a locally constant roof function over a shift space $X$.  Let $\im(r)$ be the range of $r$.  If the elements in $\im(r)$ are setwise commensurate, then the suspension flow $\ph^t$ determined by $r$ is not topologically mixing.
	\end{theorem}
	
	We show in Lemma \ref{lem:fin-range} that a locally constant function on a compact space can only take on finitely many values.  This implies the following corollary.
	
	\begin{corollary}\label{cor:rat-not-mix}
		If $r$ is a locally constant roof function over a shift space $X$, where $\im(r)\subset\QQ$, then the flow $\ph^t$ determined by $r$ is not topologically mixing.
	\end{corollary}
	
	\begin{remark}
		In \cite{mixing-bowen} and \cite{Plante}, the roof functions for a non-mixing suspension flow had to be cohomologous to a constant.  If a roof function $r$ is cohomologous to a constant $k$, then it is relatively easy to find $k$, especially with the presence of periodic points.  This is because the constant $k$ must be unique.  
		
		However, if a roof function $r$ is cohomologous to a locally constant roof function $s$, then it is much more difficult to identify such a function because there is no unique choice of $s$.  Indeed, if $r$ is locally constant, then for any locally constant transfer function $g$, the roof function $s(x) = r(x) + g(\sigma x) - g(x)$ is also locally constant.  Moreover, it is unclear whether a continuous function is always cohomologous to a locally constant one.
		
	\end{remark}
	
	\begin{remark}
		For more general shift spaces, a result like Theorem \ref{thm:SFT-cohomologous} does not hold.  In \S\ref{sec:no-mix-incommensurate}, we provide an example of a locally constant roof function over a coded subshift that produces a flow that is not topologically mixing, but the flow has periodic points with incommensurate orbit lengths. Thus, it cannot be cohomologous to a roof function that only takes values in $\delta\NN$ for any $\delta \in \RR$.
	\end{remark}
	
	Since any continuous roof function can be approximated by a locally constant roof that takes values in $\QQ$, we immediately get the following theorem.
	
	\begin{theorem}\label{thm:no-mix-dense}
		For any shift space $X$ on a finite alphabet, the set of roof functions that induce suspension flows that are not topologically mixing is dense in the set of positive continuous real-valued functions.
	\end{theorem}
	The proofs of Theorem \ref{thm:loc-const-com-nonmixing} and \ref{thm:no-mix-dense} can be found in \S\ref{sec:density}.
	
	Unlike connected topological spaces, a locally constant function defined on shift space need not be constant because shift spaces are totally disconnected.  As we see in Theorem \ref{thm:loc-const-com-nonmixing} and Corollary \ref{cor:rat-not-mix}, there is an abundance of locally constant roof functions that produce flows that are not topologically mixing but are not necessarily constant.  However, these kinds of roof functions cannot be defined on a connected metric space without losing continuity.  
	
	\begin{remark}
		It would be interesting to see if an absence of topological mixing in the flow implied that the roof function must be cohomologous to a constant for connected metric spaces, but it is currently unknown.
	\end{remark}
	
	\subsection{Open questions}\label{sec:questions}
	In Theorem \ref{thm:fish-hass-incom} and in \cite{mixing-bowen}, the periodic points are dense in the set $\Lambda$.  Although the statements in Theorems \ref{thm:mix-dich} and \ref{thm:mix-dich-beta} do not explicitly mention that the periodic points are dense, it is true in both of these settings. However, not all subshifts with dense periodic points possess a synchronizing word or are $\beta$-shifts, so there are still questions to be answered more generally.
	
	In addition to general shift spaces that have dense periodic points, it is currently unknown how to characterize a topologically mixing dichotomy for suspension flows over shift spaces without dense periodic points.  There are many examples of such shift spaces including minimal shift. Topological mixing properties of some minimal shifts were studied in \cite{DekkingKeane, Pet}.  Other examples of minimal shifts include the Morse-Thue shift, which is uniquely ergodic \cite{hB22}. There are also examples of minimal subshifts with arbitrarily many measures of maximal entropy \cite{DGS76}.    
	
	There are also topologically mixing shift spaces with only finitely many periodic points.  For instance, \cite{Kwietniak} studies a weakly topologically mixing shift space with only a single periodic orbit.  In \S\ref{sec:mix-no-per}, we have included a construction of a topologically mixing shift with only two periodic orbits.  Very little is known about the topological mixing properties of suspension flows over these kinds of shift spaces.
	
	\begin{question}
		Suppose $X$ is a transitive shift space with dense periodic points.  If $r$ is a continuous roof function where the period lengths under the flow are all setwise commensurate, is the suspension flow not topologically mixing?  Is the roof function cohomologous to a locally constant roof function?
	\end{question}
	
	\begin{question}
		Suppose $X$ is a transitive shift space with finitely many periodic points.  If $r$ is a continuous roof function that produces a suspension flow that is not topologically mixing, is it cohomologous to a constant function or a locally constant function?
	\end{question}
	
	\begin{question}
		If $X$ is a minimal shift space and $r$ is a continuous roof function that produces a suspension flow that is not topologically mixing, is it cohomologous to a constant?
	\end{question}
	
	\section{Examples}\label{sec:examples}
	\subsection{Cohomologous to a constant is too restrictive}\label{sec:no-mix-ex}
	
	Todd Fisher produced a simple example of a roof function that is not cohomologous to a constant and yields a suspension flow that is not topologically mixing.  We will refer to the suspension flow with the following roof function as $\ph^t$ in this section.
	\begin{example}\label{ex:not-mixing}
		Consider the full shift on two symbols $\Sigma_2 = \{0,1\}^\ZZ$ and define the roof function to be
		\begin{equation}
			r(x) = \begin{cases}
				2 & \text{if } x \in [0]\\
				3 &  \text{if } x \in [1].
			\end{cases}
		\end{equation}  
		The roof function $r$ is not cohomologous to a constant, and the resulting suspension flow $\ph^t$ is not topologically mixing.  
	\end{example}
	\begin{proof}
		
		Let $\overline{0}$ be the bi-infinite sequence of 0s and $\overline{1}$ be the bi-infinite sequence of 1s.  The points $\overline{0}$ and $\overline{1}$ are the fixed points of the system $(\Sigma_2,\sigma)$. Since $r(\overline{0})\neq r(\overline{1})$, it is apparent that this function is not cohomologous to a constant by Proposition \ref{prop:not-cohomologous}.
		
		\begin{figure}[!ht]
			\centering
			\def\len{2}
			\def\high{1.8}
			\def\gap{.75}
			\def\lilgap{.35}
			\def\tinygap{0}
			\begin{tikzpicture}
				\fill [gray!20] (0,0) rectangle ({\len},{\high-\tinygap});
				\fill [gray!20] (0,{\high+\tinygap}) rectangle ({\len},{2*\high});
				\fill [gray!20] ({\len+\gap},0) rectangle ({2*\len+\gap},{\high-\tinygap});
				\fill [gray!20] ({\len+\gap},{\high+\tinygap}) rectangle ({2*\len+\gap},{2*\high - \tinygap});
				\fill [gray!20] ({\len+\gap},{2*\high+\tinygap}) rectangle ({2*\len+\gap},{3*\high});
				
				\draw (0,0) -- node[below] {$[0]$} ++ ({\len},0);
				\draw ({\len+\gap},0) -- node[below] {$[1]$} ++ ({\len},0);
				
				\draw[dashed] (0,{\high}) -- ({\len},{\high});
				\draw[dashed] ({\len+\gap},{\high}) -- ({2*\len+\gap},{\high});
				\draw[dashed] ({\len+\gap},{2*\high}) -- ({2*\len+\gap},{2*\high});
				
				\draw (0,{2*\high}) -- ({\len},{2*\high});
				\draw ({\len+\gap},{3*\high}) -- ({2*\len+\gap},{3*\high});
				
				\node at ({-\lilgap},{\high}) {$1$};
				\node at ({-\lilgap},{2*\high}) {$2$};
				\node at ({-\lilgap},{3*\high}) {$3$};
				
				\path (0,{0.5*\high}) -- ({\len},{0.5*\high}) node[midway] {$A$};
				\path (0,{1.5*\high}) -- ({\len},{1.5*\high}) node[midway] {$B$};
				\path ({\len+\gap},{0.5*\high}) -- ({2*\len+\gap},{0.5*\high}) node[midway] {$C$};
				\path ({\len+\gap},{1.5*\high}) -- ({2*\len+\gap},{1.5*\high}) node[midway] {$D$};
				\path ({\len+\gap},{2.5*\high}) -- ({2*\len+\gap},{2.5*\high}) node[midway] {$E$};
			\end{tikzpicture}
			\caption{Partition of the phase space of the flow $\ph^t$.}\label{fig:susp-part}
		\end{figure}
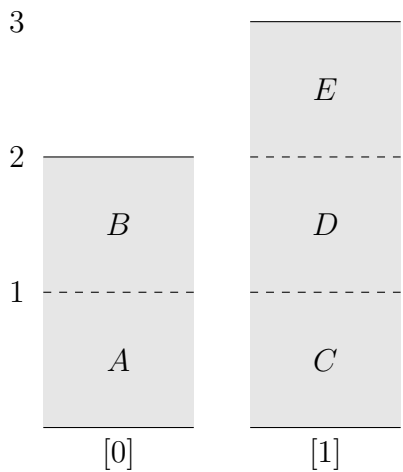
		
		We partition the phase space $M$ of $\ph^t$ into five subsets of height 1, as sketched in Figure \ref{fig:susp-part}.
		
		\begin{align*}
			A &=\{(x,y)\in M: x\in[0]\text{ and } x\in[0,1)\}\\
			B &=\{(x,y)\in M: x\in[0]\text{ and } x\in[1,2)\}\\
			C &=\{(x,y)\in M: x\in[1]\text{ and } x\in[0,1)\}\\
			D &=\{(x,y)\in M: x\in[1]\text{ and } x\in[1,2)\}\\
			E &=\{(x,y)\in M: x\in[1]\text{ and } x\in[2,3)\}\\
		\end{align*}
		
		We  associate the orbit of a point in $M$ with a coding based on which partition elements the point visits under the flow.  For example, if a point in $M$ begins in $A$ it must flow into $B$ and can either return back to $A$, or it can get mapped to the other cylinder and flow into $C$. The set of codings yields a subshift of finite type associated with the direct graph in Figure \ref{fig:SFT-graph}.
		
		\begin{figure}[h!]
			\centering
			
			\def\Len{2}
			\def\High{1.8}
			\def\Gap{.75}
			\def\Tinygap{1.5}
			\begin{tikzpicture}
				\tikzset{vertex/.style = {shape=circle,draw,minimum size=1.5em}}
				\tikzset{edge/.style = {->,> = latex'}}
				
				\node[vertex] (1) at (0,0) {$A$};
				\node[vertex] (2) at (0,{\High}) {$B$};
				\node[vertex] (3) at ({\Len+\Gap},0) {$C$};
				\node[vertex] (4) at ({\Len+\Gap},{\High}) {$D$};
				\node[vertex] (5) at ({\Len+\Gap},{2*\High}) {$E$};
				
				\draw[edge, shorten >={\Tinygap}] (1.100) to (2.260);
				\draw[edge, shorten >={\Tinygap}] (2.280) to (1.80);
				\draw[edge, shorten >={\Tinygap}] (2) to (3);
				\draw[edge, shorten >={\Tinygap}] (3) to (4);
				\draw[edge, shorten >={\Tinygap}] (4) to (5);
				\draw[edge, shorten >={\Tinygap}] (5) to (1);
				\draw[edge, shorten >={\Tinygap}, bend left = 40] (5) to (3);	
			\end{tikzpicture}
			\caption{Base dynamics of conjugate flow $\psi^t$}\label{fig:SFT-graph}
		\end{figure}
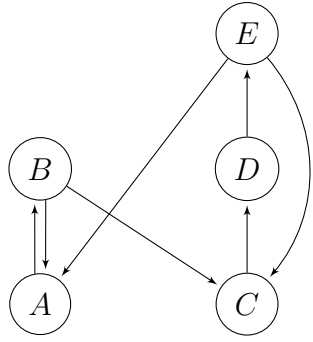
		
		It follows that $X = \{AB,CDE\}^\ZZ$ is the resulting subshift of finite type. We construct a new suspension over $X$ whose roof function is identically equal to 1 and call it $\psi^t$. The flow $\psi^t$ is conjugate to the original suspension $\ph^t$. 
		
		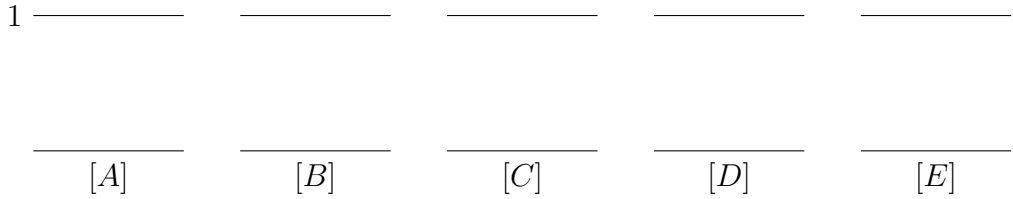
\begin{figure}[h!]
			
			\centering
			
			\def\LEN{2}
			\def\hi{1.8}
			\def\GAP{.75}
			\def\LilGap{.25}
			\begin{tikzpicture}
				\draw (0,0) -- node[below] {$[A]$} ++ ({\LEN},0);
				\draw ({\LEN+\GAP},0) -- node[below] {$[B]$} ++ ({\LEN},0);
				\draw ({2*\LEN+2*\GAP},0) -- node[below] {$[C]$} ++ ({\LEN},0);
				\draw ({3*\LEN+3*\GAP},0) -- node[below] {$[D]$} ++ ({\LEN},0);
				\draw ({4*\LEN+4*\GAP},0) -- node[below] {$[E]$} ++ ({\LEN},0);
				
				\draw (0,{\hi}) -- node[below]{} ++ ({\LEN},0);
				\draw ({\LEN+\GAP},{\hi})-- node[below]{} ++ ({\LEN},0);
				\draw ({2*\LEN+2*\GAP},{\hi}) -- node[below]{} ++ ({\LEN},0);
				\draw ({3*\LEN+3*\GAP},{\hi}) -- node[below]{} ++ ({\LEN},0);
				\draw ({4*\LEN+4*\GAP},{\hi}) -- node[below]{} ++ ({\LEN},0);
				
				\node at ({-\LilGap},{\hi}) {$1$};
			\end{tikzpicture}
			\caption{$\psi^t$ is a suspension flow with roof function $1$.}
		\end{figure}
		
		It is well known that a suspension flow with a constant roof function is not topologically mixing.  Since $\ph^t$ is conjugate to $\psi^t$, $\ph^t$ is not topologically mixing.  
	\end{proof}
	
	This counterexample shows that a suspension flow can have a roof function that is not cohomologous to a constant, but also fails to be topologically mixing.
	
	An alternative way to think about this example is that we are constructing a new suspension flow by taking the appropriate transverse cross-section of the original flow $\ph^t$.  If we choose $\Sigma_2 \times \{0\} \cup \Sigma_2 \times \{1\} \cup [1]\times\{2\}$ as our transverse cross-sections, then we would obtain $\psi^t$ because the return time to this cross-section is always 1, and the return map defined on the cross-section would be governed by $X$.
	
	\subsection{Pairwise commensurate periodic orbits is insufficient}\label{sec:mixing-example}
	Suppose $P$ is a finite subset of $\RR$. If $x/y \in \QQ$ for any pair $x,y\in P$, then we know that $P$ is a setwise commensurate set.  However, if $P$ is countably infinite, then $x/y \in \QQ$ for any pair $x,y\in P$, may not imply that $P$ is a setwise commensurate set.  The following example demonstrates the importance of this distinction in our context.
	
	\begin{example}\label{ex:mixing}
		Let $\Sigma_2 = \{0,1\}^\ZZ$ and define $\rho\colon \Sigma_2 \to \ZZ \cup\{\infty\}$ to be the number of consecutive 0s starting from the 0th position.  That is, $\rho(\dots.0^n1\dots) = n$.  Define $r\colon \Sigma_2 \to \RR$ by 
		\[
		r(x) = 
		\begin{cases}
			1 + \frac{1}{1+\rho(x)} & \text{if } x\in [0]\\
			1 & \text{if } x \in [1]
		\end{cases}
		\] 
		and note that $\rho(\overline{0}) = 1$.  The roof function $r$ is continuous and any pair of periodic points have commensurate orbit lengths (their ratio is rational); however, the set of all period lengths is not setwise commensurate in the sense of Definition \ref{def:commensurate}.  The suspension flow with roof function $r$ is topologically mixing.
	\end{example}
	\begin{proof}
		Let $u$ and $v$ be arbitrary words of equal length, $|u|=|v|$.  Without loss of generality, we may consider $U = [u] \times (0,\eps)$ and $V = [v] \times (0,\eps)$ to prove topological mixing by Lemma \ref{lem:mixing-simple}.  Consider the collection of points 
		\[x_{m,n} = \dots u 1 0^m 1^n v^\infty \in [u].\]
		Note that we have not indicated the center of $x_{m,n}$ using ``." in the usual way.  However, since we are requiring that $x_{m,n} \in [u]$, it is understood that the ``." would be placed in the middle of $u$.  Also note that for any $n,m$ and $s\in(0,\eps)$ the point $(x_{m,n},s)$ starts in $U$, eventually makes its way to $V$, and then periodically returns to $V$ after some uniform time.  We denote this return time as $\omega(v)$ and note that $\omega(v)$ is independent of $n$ and $m$ since $r$ is locally constant away from $\overline{0}$ and $v$ has finite length.
		
		In order to show that the flow is topologically mixing, we must show that the time it takes $x_{m,n} \times (0,\eps)$  to enter $V$ mod $\omega(v)$ for all $m,n$ is dense in $[0,\omega(v)]$ by Lemma \ref{lem:mixing-simple}.  
		
		Note that $\sum_{j=0}^{|u|}r(\sigma^j x_{m,n})$ is constant for any pair of $m,n$.  We denote this constant by $\omega(u1)$.  This means that the time it takes for the line segment $x_{m,n}\times(0,\eps)$ to be entirely contained in $V$ is 
		\begin{equation}\label{eqn:tau-mn}
			\begin{split}
				\tau(x_{m,n}) & \coloneqq  \sum_{j=0}^{|u|+m+n}r(\sigma^j x_{m,n}) 
				= \omega(u1) + \sum_{j=1}^m\left(1+\frac{1}{1+j}\right) + n\\
				& = \omega(u1) + m + n +\sum_{j=2}^{m+1}\frac{1}{j}.
			\end{split}
		\end{equation}
		That is, $\ph^{\tau(x_{m,n})}(x_{m,n}\times(0,\eps)) \cap V = \ph^{\tau(x_{m,n})}(x_{m,n}\times(0,\eps))$.  Note that here we are using the fact that $|u| = |v|$.  If $|u| \neq |v|$, then $\sigma^{|u|+1+m+n}x_{m,n}$ may not be contained in $[v]$.  Additionally, we have that for all $N \in \NN$
		\[
		\ph^{\tau(x_{m,n}) + N\omega(v)}(x_{m,n}\times(0,\eps)) \cap V = \ph^{\tau(x_{m,n}) + N\omega(v)}(x_{m,n}\times(0,\eps)).
		\]
		This justifies the reasoning that we only need to show that the collection of $\tau(x_{m,n})$ mod $\omega(v)$ is dense in $[0,\omega(v)]$.
		
		Let $\delta > 0$ and $\alpha \in (0,\omega(v))$.  There exists an $m$ such that 
		\[
		\left|\omega(u1) + \sum_{j=2}^{m+1} \frac{1}{j} - \alpha \mod \omega(v)\right| < \delta
		\]
		since the sum diverges but eventually increases by arbitrarily small increments.  
		
		Note that $\omega(v)$ is a rational number of the form $p/q$ where $p$ and $q$ are coprime.  Hence, there exists an $n$ (which depends on $m$) so that 
		\[
		\frac{n+m}{p/q} \in \NN \iff n + m \equiv 0 \mod \omega(v).
		\]
		Thus by \eqref{eqn:tau-mn} we have 
		\begin{align*}
			|\tau(x_{m,n}) - \alpha \mod \omega(v)| & = \left| \omega(u1) + m + n +  \sum_{j=2}^{m+1} \frac{1}{j} - \alpha \mod \omega(v) \right|\\
			& = \left|\omega(u1) + \sum_{j=2}^{m+1} \frac{1}{j} - \alpha \mod \omega(v) \right|
			< \delta.
		\end{align*}
		
		Since $\delta$ is arbitrary it follows that the collection of $\tau(x_{m,n})$ is dense in $[0,\omega(v)]$.
		
		Observe that once $x_{m,n}\times(0,\eps)$ has flowed into $V$, it returns to $V$ after time $\omega(v)$.  It also intersects $V$ for a duration of $2\eps$ time.  We can cover the interval $[0,\omega(v)]$ with finitely many intervals of length $2\eps$ whose centers correspond to $\tau(x_{m,n})$ because the collection of $\tau(x_{m,n})$ is dense. Therefore, the suspension is topologically mixing.
	\end{proof}
	\begin{remark}
		The approach used to prove topological mixing here is a simplified version of the proof of topologically mixing for more general roof functions that appear in \S\ref{sec:mix-dict-pf} and \S\ref{sec:beta-shifts}.
	\end{remark}
	
	\subsection{A non-mixing flow with incommensurate periodic orbit lengths}\label{sec:no-mix-incommensurate}
	
	\begin{definition}
		Let $C$ be a set of words from an alphabet $A$.  Define 
		\[
		B_C = \{x \in A^{\ZZ} : \exists \text{ a bi-infinite increasing sequence } s_n \text{ s.t. } \forall k, \, x_{[s_n,s_{n+1})} \in C\},
		\]
		and let $X_C = \overline{B_C}$.  The set $X_C$ is the \emph{coded subshift} associated to $C$.
	\end{definition}
	Heuristically, a coded subshift is the closure of all arbitrary bi-infinite concatenations of words from the generating set $C$.  They were first introduced in \cite{BH86} and are a broad family that includes shift spaces like $S$-gap shifts and $\beta$-shifts.  
	
	\begin{example}\label{ex:no-mix-incommensuarte}
		Let $C = \{2^n3^n : n \in \NN\} \cup \{0,1\}$ and $X_C \subset \{0,1,2,3\}^\ZZ$ be the coded subshift associated to $C$.  Let $a , b > 0$ be a pair of real numbers such that $a$ and $b$ are incommensurate; that is, $a/b \notin\QQ$.  Define
		\[
		r(x) = \begin{cases}
			a + b & \text{if } x\in [0] \cup [1] \\
			a & \text{if } x \in [2] \\
			b & \text{if } x \in [3]
		\end{cases}
		\]
		The suspension flow with roof $r$ is not topologically mixing, but there does not exist a $\delta \in \RR$ satisfying equation \eqref{eqn:no-mix-orbit-length} for all periodic points in $X_C$.
	\end{example}
	\begin{proof}
		Note that $X_C$ is a topologically mixing shift with synchronizing words, but it is not a subshift of finite type.  Observe that $\overline{2}$ and $\overline{3}$ are both elements of $X_C$.  If $x \neq \overline{2}, \overline{3}$ and is a periodic point, then 
		\[\sum_{j=0}^{\per(x)-1}r(\sigma^j x) \in (a+b)\ZZ.\]  
		Furthermore, the symbol $0$ is a synchronizing word and $\overline{2},\overline{3} \notin [0]$, so the suspension flow is not topologically mixing by Theorem \ref{thm:mix-dich}.
		
		However, $r(\overline{2})$ and $r(\overline{3})$ are incommensurate, so we cannot find a $\delta$ that satisfies the setwise commensurate condition in equation \eqref{eqn:no-mix-orbit-length} for all the periodic points in $X_C$.
	\end{proof}
	
	Theorem \ref{thm:SFT-cohomologous} asserts that if a suspension flow over a subshift of finite type is not topologically mixing, then there exists a cross-section transverse to the flow such that the return time to this cross-section is constant, which is evident in Example \ref{ex:not-mixing}.  However, there is no such cross-section for the suspension flow in Example \ref{ex:no-mix-incommensuarte}.
	
	\begin{remark}
		Note that Example \ref{ex:no-mix-incommensuarte} does not contradict Theorem \ref{thm:mix-dich}.  Although cylinders of the form $[2^n]$ possess pairs of periodic points with incommensurate orbit lengths, $2^n$ is not a synchronizing word for any $n \in \NN$.  For instance, $02^n$ and $2^n3^{n+1}0$ are permissible words, but $02^n3^{n+1}0$ is a forbidden word.  Similar for words of the form $3^n$.   This motivates the following proposition.
	\end{remark}
	
	\begin{proposition}\label{prop:mult-sync}
		Suppose $X$ is a transitive shift space with a synchronizing word $v$ and $r \colon X \to (0,\infty)$ is a roof function satisfying the Walters property. Additionally, suppose there exists a $\delta > 0$ such that for any periodic point $p \in [v]$, we have
		\[
		\sum_{j=0}^{\per(p)-1} r(\sigma^j p) \in \delta \ZZ.
		\]
		If $u$ is another synchronizing word, then for any periodic point $q \in [u]$
		\[
		\sum_{j=0}^{\per(q)-1} r(\sigma^j q) \in \delta \ZZ.
		\]
	\end{proposition}
	\begin{proof}
		Since $q \in [u]$ is periodic there exists a word $w = uw' $ such that $q = \overline{w}$ (note that $w'$ could be the empty word).  Since $X$ is transitive, there are words $v'$ and $v''$ such that $vv'u \in \LLL$ and $uw'v''v \in \LLL$.  Since $u$ is synchronizing and $w^n \in \LLL$ for all $n \in \NN$, it follows that $vv'w^n \in \LLL$ and $w^n v''v \in \LLL$ for all $n\in\NN$.  Thus, $vv'w^nv''v \in \LLL$, for all $n \in \NN$ as well.
		
		Let $\eps > 0$ and $k$ be the constant required to use the Walters property.  There exist periodic points $x,y \in [v]$ of the form
		\[
		x  = \overline{vv'w^{2k+1}v''} \qquad\text{and}\qquad y  = \overline{vv'w^{2k}v''}.
		\]
		
		Since $x$ and $y$ are periodic points in $[v]$, there exists $K_1, K_2 \in \NN$ such that
		\[
		K_1 \delta = \sum_{j=0}^{\per(x)-1} r(\sigma^j x) \qquad\text{and}\qquad
		K_2 \delta =  \sum_{j=0}^{\per(y)-1} r(\sigma^j y).
		\]
		Applying the Walters property to different parts of the orbit of $x$ we get
		\begin{align}
			\left|\sum_{j=0}^{|vv'w^k|-1} r(\sigma^j x) - \sum_{j=0}^{|vv'w^k|-1} r(\sigma^j y) \right| & < \eps, \label{eqn:x-y-1} \\ 
			\left|\sum_{j=|vv'w^k|}^{|vv'w^{k+1}|-1} r(\sigma^{j} x) - \sum_{j=0}^{|w|-1} r(\sigma^j q)\right| & < \eps, \label{eqn:x-z} \\
			\left|\sum_{j=|vv'w^{k+1}|}^{|vv'w^{2k+1}v''|-1} r(\sigma^{j} x) - \sum_{j=|vv'w^{k}|}^{|vv'w^{2k}v''|-1} r(\sigma^{j} y)\right| & < \eps . \label{eqn:x-y-2}
		\end{align}
		Using equations \eqref{eqn:x-y-1} and \eqref{eqn:x-y-2} we have
		\begin{equation}\label{eqn:x-y-3}
			\left|\sum_{j=0}^{|vv'w^k|-1} r(\sigma^j x) + \sum_{j=|vv'w^{k+1}|}^{|vv'w^{2k+1}v''|-1} r(\sigma^{j} x) - \sum_{j=0}^{\per(y)-1} r(\sigma^{j} y)\right| < 2\eps.
		\end{equation}
		Using equations \eqref{eqn:x-z} with \eqref{eqn:x-y-3} we get
		\begin{align*}
			\left|(K_1 - K_2)\delta - \sum_{j=0}^{\per(q)-1} r(\sigma^j q)\right| & = \left|\sum_{j=0}^{\per(x)-1} r(\sigma^j x) - \sum_{j=0}^{\per(y)-1} r(\sigma^j y) - \sum_{j=0}^{\per(q)-1} r(\sigma^j q)  \right| \\
			& < \left|\sum_{j=|vv'w^k|}^{|vv'w^{k+1}|-1} r(\sigma^{j} x) - \sum_{j=0}^{|w|-1} r(\sigma^j q)\right| + 2\eps \\
			& < 3\eps.
		\end{align*}
		Since $\eps$ was arbitrary, we get that
		\[
		\sum_{j=0}^{\per(q)-1} r(\sigma^j q) \in \delta \ZZ.
		\]
	\end{proof}
	
	Proposition \ref{prop:mult-sync} ensures that a roof function cannot satisfy the setwise commensurate condition in equation \eqref{eqn:no-mix-orbit-length} on the cylinder of one synchronizing word, but violate it on the cylinder of some other synchronizing word. 
	
	This implies that one cannot obtain the phenomenon in Example \ref{ex:no-mix-incommensuarte} for subshifts of finite type because all sufficiently long words of a subshift of finite type are synchronizing.
	
	\begin{corollary}\label{cor:SFT-commensurate}
		If a suspension flow over a subshift of finite type is not topologically mixing, then the set of all periodic orbit lengths of the flow is setwise commensurate.
	\end{corollary}
	
	\section{Proof of Theorem \ref{thm:mix-dich}}\label{sec:mix-dict-pf}
	
	It is enough to consider sets of the form $[w] \times (a,b)$ to verify that the suspension flow is topological mixing.  These sets are open under the topology of the phase space $M_r$ of the suspension.  Although not all sets are of this form, every open set contains a subset of the form $[w]\times (a,b)$.  Thus, it is sufficient to prove topological mixing using sets like these.  
	
	We will make an additional simplification to these sets.  We will consider sets of the form $[w] \times (0,\eps)$, for $\eps>0$.  This is mainly for convenience, as topological mixing on sets like this implies topological mixing on sets $[w] \times (a, a+\eps)$ by simply adjusting the times.

	\begin{lemma}\label{lem:mixing-simple}
		Let $u,w \in \LLL$ and $0 < \eps < \min r$. If there exists a $T$ such that for all $t \geq T$, $\ph^t([u]\times (0,\eps)) \cap ([w]\times (0,\eps)) \neq \emptyset$, then for any $0<a < \min r|_{[u]}-\eps$ and $0<b < \min r|_{[w]}-\eps$, there exists a $T'$ such that for all $t\geq T'$ we have $\ph^t([u]\times (a,a+\eps)) \cap [w]\times (b,b+\eps) \neq \emptyset$.
	\end{lemma}
	\begin{proof}
		For convenience, let $U = [u]\times (0,\eps)$, $W = [w]\times (0,\eps)$, $U' = [u]\times (a,a+\eps)$, and $W' = [w]\times (b,b+\eps)$.  For $t \geq T$, there exists a point $x \in U$ such that $\ph^t(x) \in W$.  Note that here $x \notin [u]$ as it is not in the shift space.  Rather, it is a point in the actual phase space of the flow. 
		
		Observe that $\ph^a(x) \in U'$ and $\ph^{t+b}(x) \in W'$.  This implies that $\ph^{t+b-a}(\ph^{a}(x)) \in W'$, so for all $s \geq T +b -a$, we have that $\ph^s(U') \cap W' \neq \emptyset$.
	\end{proof}
	
	To prove the sufficient condition of Theorem \ref{thm:mix-dich} we need the following lemma.  Here we will let $\gcd(x,y)$ be the largest real number $c$ such that $x/c,y/c \in \NN$.  If no such number exists, then we set $\gcd(x,y)=0$.
	\begin{lemma}\label{lem:gcd-dense}
		Let $a,b\in\RR$ with $a,b\geq 0$ and suppose that their  $\gcd(a,b)<\delta$.  If $\omega > \delta$ and $x \in [0,\omega)$, then there are $n,m\in\NN$ such that $|na+mb - x \mod \omega|< \delta$.  Moreover, $n,m$ can be chosen to be arbitrarily large.
	\end{lemma}
	\begin{proof}
		Without loss of generality, we will assume that $a,b < \omega$.  There are two cases: either $a$ or $b$ is incommensurate with $\omega$, or both $a$ and $b$ are commensurate with $\omega$.
		
		If $a$ is incommensurate with $\omega$, then the set $\{na \mod \omega: n\in\NN\}$ is dense in $[0,\omega)$ as this action is equivalent to an irrational rotation of a circle of circumference $\omega$.  Thus for any $m \in \NN$ it follows that $\{na +mb \mod \omega : n \in \NN\}$ is dense in $[0,\omega)$.  Similarly when $b$ is incommensurate with $\omega$.
		
		Now consider the case when $a$, $b$, $\omega$ are commensurate, with gcd $0 <  \delta' \leq \delta$ where $\delta$ is still the $\gcd(a,b)$.  Since $a/\delta'$, $b/\delta'$, and $\omega/\delta'$ are all natural numbers, it follows that there exist $n,m \in \NN$ such that
		\[\frac{na}{\delta'} + \frac{mb}{\delta'} \equiv \frac{\delta}{\delta'} \mod \frac{\omega}{\delta'},\]
		since $\delta/\delta'$ is the gcd of $a/\delta'$ and $b/\delta'$.  This is equivalent to saying that there exists an $N \in \NN$ such that
		\[\frac{na}{\delta'} + \frac{mb}{\delta'} - \frac{N\omega}{\delta'} = \frac{\delta}{\delta'}.\]
		Multiplying by $\delta'$, we get
		\[na + mb - N\omega = \delta.\]
		So for any $x \in [0,\omega)$ we can find a $k \in \NN$ such that
		\[|kna + kmb - x -kN\omega| < \delta,\]
		which completes the proof.
	\end{proof}
	
	The primary consequence of Lemma \ref{lem:gcd-dense} is that the set $\{na+mb \mod \omega : n,m \in \NN\}$  is dense in $[0,\omega)$, or neighboring points are at most distance $\delta$ apart.  We will apply this lemma when the lengths of periodic points under the flow is setwise incommensurate.  If a set is setwise incommensurate, then it generates a dense subset of $\RR$ by definition; however, it is important in our setting that the coefficients are positive.
	
	We are now ready to prove Theorem \ref{thm:mix-dich}.
	
	\begin{proof}
		We begin by proving the necessary condition.  That is, if $r$ satisfies the setwise commensurate condition in equation \eqref{eqn:no-mix-orbit-length}, then the suspension flow is not topologically mixing.  Let $r$ be a roof function satisfying the Walters property over a subshift $X$ with a synchronizing word $v$.  By hypothesis, there exists a $\delta > 0$ such that
		\begin{equation}
			\sum_{j=0}^{\per(x)-1} r(\sigma^j x) \in \delta\ZZ,
		\end{equation}	
		
		Without loss of generality, we assume that $\delta \ll \min r$.  Let $0 < \eps < \delta / 8$, and let $k$ be the length needed to satisfy the Walters property for this value of $\eps$.  There is a word $u$ such that $vu \in \LLL$ and $|vu|> 2k$.  Consider the set $V = [vu] \times (0,\eps)$.
		
		Suppose $x \in [vu]$ and $\sigma^n x \in [vu]$.  Then there exists a $w\in\LLL$ such that $|vuw| = n$ and $x = \dots vuwvuw'\dots$.  Here $w'$ is not any special word in particular; it is meant to signify that after $vuwvu$, we do not know what the tail of $x$ looks like.  
		
		Since $vuwv \in \LLL$ and $v$ is synchronizing, it follows that $\overline{vuw} \in X$.  Moreover,
		\begin{equation}\label{eqn:return-near-delta-scale}
			\left| \sum_{j=0}^{n-1} r(\sigma^j x) - r(\sigma^j \overline{vuw}) \right| < \eps
		\end{equation}
		by the Walters property.  Since $\overline{vuw}$ is periodic, its orbit has length $m\delta$ for some $m \in \NN$.  Let 
		\[\tau_x = \sum_{j=0}^{n-1} r(\sigma^j x).\]
		
		Since $\ph^t(x\times(0,\eps)) \cap V \neq\emptyset$ for $t \in (\tau_x-\eps, \tau_x+\eps)$, by equation \eqref{eqn:return-near-delta-scale}, it follows that $m\delta \in (\tau_x-\eps,\tau_x+\eps)$.  In other words, $(\tau_x-\eps,\tau_x+\eps) \cap \delta\ZZ \neq\emptyset$.
		
		Since $x$ and the return of $x$ to $[vu]$ was arbitrary, it follows that for any $x \in [vu]$, the return of $x \times (0,\eps)$ to $V$ under the flow must be during an interval of length $2\eps$ intersecting $\delta\ZZ$.  Note that by our choice of $\eps$ we know that $2\eps < \delta/4$.  Hence, the collection of times $i\delta + \delta/2$ for $i \in \NN$ cannot be included in any of these return intervals.  Therefore, there is no such $T$ such that $\ph^t(V) \cap V \neq \emptyset$ for all $t > T$, and the flow is not topologically mixing.
		
		We will now prove the sufficient condition of Theorem \ref{thm:mix-dich} and do this by contraposition.  That is, we will show that if there is no $\delta$ that satisfies equation \eqref{eqn:no-mix-orbit-length}, then the suspension flow must be topologically mixing.  The proof is essentially an adaptation of the proof of topological mixing in Example \ref{ex:mixing}.
		
		Let $u,w \in \LLL$ be any pair of words and $\eps > 0$ of equal length.  Let $U = [u] \times (0,\eps)$ and $W = [w]\times (0,\eps)$.  Fix $0 < \gamma < \eps/21$.  There exists a $k$ such that if $x,y\in X$ satisfying $x_{[-k,n+k]}=y_{[-k,n+k]}$, then $|\sum_{j=0}^{n-1}r(\sigma^j x) - \sum_{j=0}^{n-1}r(\sigma^j y)|<\gamma$ by the Walters property.
		
		As before, let $v$ be the synchronizing word of $X$.  By transitivity there are words $w'$ and $u'$ such that $ww'v\in\LLL$ and $vu'u\in\LLL$.  
		
		For any periodic point $p$, let $\omega(p)$ denote the period of $(p,0)$ under the suspension flow $\ph^t$.  Let us assume that condition \eqref{eqn:no-mix-orbit-length} fails, that is, by hypothesis, there exist periodic points $p,q \in [v]$ such that either $\omega(p)$ and $\omega(q)$ are incommensurate or $\gcd(\omega(p),\omega(q)) \leq \gamma$.  Consequently, there exists $\alpha,\beta \in \NN$ such that 
		\[0 < |\alpha\omega(p)-\beta\omega(q)| \leq \gamma.\]
		
		Since $p,q \in [v]$ are periodic, there exists words $v_1,v_2 \in \LLL$ such that $p = \overline{v_1} = \overline{vv'}$ and $q = \overline{v_2} = \overline{vv''}$.  Since $vv'v,vv''v \in \LLL$, it follows that $v_1^nv_2^m \in \LLL$ for any $n,m \in \NN$.  Note that here we are slightly abusing notation as $v_1$ and $v_2$ are words and do not represent symbols of $v$.
		
		Additionally, there is a word $u''$ such that $vu'uu''v \in \LLL$ by transitivity.  Hence $vu'uu'' vu'uu'' \dots vu'uu'' \in \LLL$ because $v$ is synchronizing.  Let $\zeta = vu'uu''$, and consider the family of points
		\[
		x_{n,m} = \dots w^* w w' v_1^n v_2^m \zeta^\infty \in [w]
		\]
		for some permissible $w^*\in\LLL_{\geq k}$ that may precede $ww'v$.  Note that we are omitting the ``." designating the center of $x_{n,m}$ for convenience since it is understood that $x_{n,m} \in [w]$.
		
		Fix $n_0 = m_0 = 3k$ and consider $x^* \coloneqq x_{n_0,m_0}$.  Now let $n,m \geq 3k$.  By shadowing the orbit of $x_{n,m}$ with the appropriate points, we can repeatedly apply the Walters property as follows.
		\begin{equation}\label{eqn:first-est}
			\left|\sum_{j=0}^{|ww'|+k|v_1|-1} r(\sigma^j x_{n,m}) - \sum_{j=0}^{|ww'|+k|v_1|-1} r(\sigma^j x^*) \right| < \gamma
		\end{equation}
		Let $x_{n,m}^{(1)} = \sigma^{|ww'|+k|v_1|}x_{n,m}$.
		\begin{equation}\label{eqn:shadow-p}
			\left|\sum_{j=0}^{(n-2k)|v_1|-1} r(\sigma^{j} x_{n,m}^{(1)}) - \sum_{j=0}^{(n-2k)|v_1|-1} r(\sigma^j p) \right| < \gamma
		\end{equation}
		Let $x_{n,m}^{(2)} = \sigma^{(n-2k)|v_1|}x_{n,m}^{(1)}$.
		\begin{equation}\label{eqn:shadow-second}
			\left|\sum_{j=0}^{k(|v_1|+|v_2|)-1} r(\sigma^{j} x_{n,m}^{(2)}) - \sum_{j=0}^{k(|v_1|+|v_2|)-1} r(\sigma^{j+|ww'|+2k|v_1|} x^*) \right| < \gamma
		\end{equation}
		Let $x_{n,m}^{(3)} = \sigma^{k(|v_1|+|v_2|)}x_{n,m}^{(2)}$.
		\begin{equation}\label{eqn:shadow-q}
			\left|\sum_{j=0}^{(m-2k)|v_2|-1} r(\sigma^{j} x_{n,m}^{(3)}) - \sum_{j=0}^{(m-2k)|v_2|-1} r(\sigma^{j} q) \right| < \gamma
		\end{equation}
		Let $x_{n,m}^{(4)} = \sigma^{(m-2k)|v_2|}x_{n,m}^{(3)}$.
		\begin{equation}\label{eqn:arrive-est}
			\left|\sum_{j=0}^{k(|v_2|+|\zeta|)-1} r(\sigma^{j} x_{n,m}^{(4)}) - \sum_{j=0}^{k(|v_2|+|\zeta|)-1} r(\sigma^{j+|ww'|+3k|v_1|+2k|v_2|} x^*) \right| < \gamma
		\end{equation}
		Let $x_{n,m}^{(5)} = \sigma^{k(|v_2|+|\zeta|)}x_{n,m}^{(4)}$ and $N\in\NN$.
		\begin{equation}\label{eqn:return-to-zeta-time}
			\left|\sum_{j=0}^{N-1} r(\sigma^{j} x_{n,m}^{(5)}) - \sum_{j=0}^{N-1} r(\sigma^{j} \overline{\zeta}) \right| < \gamma.
		\end{equation}
		\begin{figure}[h!]
	\def\Space{1.2}
	\def\Lilspace{.1}
	\def\vertgap{0.5}
	\centering
	\begin{tikzpicture}
		\draw ({2*\Space},0) -- node[above] {$w^*ww'$} ++({\Space-\Lilspace},0);
		\draw ({3*\Space},0) -- node[above] {$v_1^k$} ++({\Space-\Lilspace},0);
		\draw ({4*\Space},0) -- node[above] {$v_1^{n-2k}$} ++({1.5*\Space-\Lilspace},0);
		\draw ({5.5*\Space},0) -- node[above] {$v_1^k$} ++({\Space-\Lilspace},0);
		\draw ({6.5*\Space},0) -- node[above] {$v_2^k$} ++({\Space-\Lilspace},0);
		\draw ({7.5*\Space},0) -- node[above] {$v_2^{m-2k}$} ++({1.5*\Space-\Lilspace},0);
		\draw ({9*\Space},0) -- node[above] {$v_2^k$} ++({\Space-\Lilspace},0);
		\draw ({10*\Space},0) -- node[above] {$(vu'uu'')^k$} ++({2*\Space-\Lilspace},0);
		\draw ({12*\Space},0) -- node[above] {$\zeta$} ++({\Space-\Lilspace},0);
		
		\node [right=0.5pt of {({13*\Space},0)}] {$\cdots$};
		\node [right=0.5pt of {({13*\Space},{-\vertgap})}] {$\cdots$};
		
		\draw ({2*\Space},{-\vertgap}) -- node[below] {$x^*$} ++({2*\Space-\Lilspace},0);
		\draw ({4*\Space},{-\vertgap}) -- node[below] {$p$} ++({1.5*\Space-\Lilspace},0);
		\draw ({5.5*\Space},{-\vertgap}) -- node[below] {$\sigma^{|ww'v_1^{2k}|}x^*$} ++({2*\Space-\Lilspace},0);
		\draw ({7.5*\Space},{-\vertgap}) -- node[below] {$q$} ++({1.5*\Space-\Lilspace},0);
		\draw ({9*\Space},{-\vertgap}) -- node[below] {$\sigma^{|ww'v_1^{3k}v_2^{2k}|}x^*$} ++({3*\Space-\Lilspace},0);	
		\draw ({12*\Space},{-\vertgap}) -- node[below] {$\overline{\zeta}$} ++({\Space-\Lilspace},0);
		
		\draw[dashed] ({4*\Space-0.5*\Lilspace},0) -- ({4*\Space-0.5*\Lilspace},{-2*\vertgap});
		\draw[dashed] ({5.5*\Space-0.5*\Lilspace},0) -- ({5.5*\Space-0.5*\Lilspace},{-2*\vertgap});
		\draw[dashed] ({7.5*\Space-0.5*\Lilspace},0) -- ({7.5*\Space-0.5*\Lilspace},{-2*\vertgap});
		\draw[dashed] ({9*\Space-0.5*\Lilspace},0) -- ({9*\Space-0.5*\Lilspace},{-2*\vertgap});
		\draw[dashed] ({12*\Space-0.5*\Lilspace},0) -- ({12*\Space-0.5*\Lilspace},{-2*\vertgap});
		
		\path ({2*\Space},{-3*\vertgap}) -- ({4*\Space-\Lilspace},{-3*\vertgap}) node[midway] {\eqref{eqn:first-est}};
		\path ({4*\Space},{-3*\vertgap}) -- ({5.5*\Space-\Lilspace},{-3*\vertgap}) node[midway] {\eqref{eqn:shadow-p}};
		\path ({5.5*\Space},{-3*\vertgap}) -- ({7.5*\Space-\Lilspace},{-3*\vertgap}) node[midway] {\eqref{eqn:shadow-second}};
		\path ({7.5*\Space},{-3*\vertgap}) -- ({9*\Space-\Lilspace},{-3*\vertgap}) node[midway] {\eqref{eqn:shadow-q}};
		\path ({9*\Space},{-3*\vertgap}) -- ({12*\Space-\Lilspace},{-3*\vertgap}) node[midway] {\eqref{eqn:arrive-est}};
		\path ({12*\Space},{-3*\vertgap}) -- ({13*\Space-\Lilspace},{-3*\vertgap}) node[midway] {\eqref{eqn:return-to-zeta-time}};
	\end{tikzpicture}
	\caption{Breakdown of how the orbit of $x_{n,m}$ is shadowed by the various points accompanied by the equation that each orbit segment corresponds to.  The upper line represents the orbit of $x_{n,m}$, and the lower line represents the orbit segments of the shadowing points.}
\end{figure}		

		Recall, $\omega(\overline{\zeta}) = \sum_{j=0}^{|\zeta|-1}r(\sigma^j\overline{\zeta})$.  By equation \eqref{eqn:return-to-zeta-time}, for any $n,m\geq 3k$ we have
		\begin{equation}\label{eqn:return-to-zeta-N}
			\left| \sum_{j=0}^{N|\zeta|-1}r(\sigma^jx_{n,m}^{(5)}) -N\omega(\overline{\zeta}) \right| < \gamma
		\end{equation}
		This implies that all of the subsequent return times of the orbit $x_{n,m}^{(5)}$ to $U$ is in the interval $(\omega(\overline{\zeta}) - \gamma, \omega(\overline{\zeta}) + \gamma)$.

		Let 
		\begin{align*}
			K & =  \sum_{j=0}^{|ww'|+k|v_1|-1} r(\sigma^j x^*) + \sum_{j=0}^{k(|v_1|+|v_2|)-1} r(\sigma^{j+|ww'|+2k|v_1|} x^*)\\
			 & \qquad + \sum_{j=0}^{k(|v_2|+|\zeta|)-1} r(\sigma^{j+|ww'|+3k|v_1|+2k|v_2|} x^*).
		\end{align*}
		Note that $K$ is independent of $n,m$ as $x^*$ is fixed.  Moreover, $K$ represents the sum of the orbit lengths of $x^*$ that are used to shadow $x_{n,m}$.  Note that they appear in equations \eqref{eqn:first-est}, \eqref{eqn:shadow-second}, and \eqref{eqn:arrive-est}.
		
		Let $T_{n,m} = |ww'v_1^nv_2^m \zeta^k|$.  If we use equations \eqref{eqn:first-est}-\eqref{eqn:arrive-est} and the triangle inequality, we have
		\begin{equation}\label{eqn:full-arrive-est}
			\left| \sum_{j=0}^{T_{n,m}-1}r(\sigma^jx_{n,m}) - (K + (n-2k)\omega(p) + (m-2k)\omega(q))\right| < 5\gamma
		\end{equation}
		for any pair $n,m$.  
		
		Let $\tau \in (0,\omega(\overline{\zeta}))$ be arbitrary.  By Lemma \ref{lem:gcd-dense}, we know that there exists $n,m,C\in\NN$ such that 
		
		\[
		\left| \tau + C\omega(\overline{\zeta}) - (K + (n-2k)\omega(p) + (m-2k)\omega(q)) \right|  < \gamma.
		\]
		Applying equation \eqref{eqn:full-arrive-est} we have
		\[\left| \tau + C\omega(\overline{\zeta}) - \sum_{j=0}^{T_{n,m}-1} r(\sigma^jx_{n,m}) \right| < 6\gamma. \]
		and \eqref{eqn:return-to-zeta-N} estimates the subsequent returns to $U$ by
		\[\left| \tau + (C+N)\omega(\overline{\zeta}) - \sum_{j=0}^{T_{n,m}+N|\zeta|-1}r(\sigma^jx_{n,m}) \right|  < 7\gamma.\]
		
		For all $N\in\NN$, the time $\sum_{j=0}^{T_{n,m}+N|\zeta|-1} r(\sigma^jx_{n,m})$ corresponds to the moment that the image of $x_{n,m}\times(0,\eps) \subset W$ under the flow is contained in $U$.  Note that here we are using the assumption that $|u|=|w|$.  Because $\gamma < \eps/21$ we are guaranteed that $x_{n,m}\times(0,\eps)$ is passing through $U$ from time $\tau + (C+N)\omega(\overline{\zeta})-\eps/3$ to $\tau + (C+N)\omega(\overline{\zeta}) + \eps/3$ for all $N \in \NN$.
		
		In other words, if we let $T = \sum_{j=0}^{T_{n,m}-1}r(\sigma^jx_{n,m})$, then for all  $t \in (\tau + (C+N)\omega(\overline\zeta) - \eps/3, \tau + (C+N)\omega(\overline\zeta) + \eps/3)$ we are guaranteed that $\ph^{T+t}(W) \cap U \neq \emptyset$ for all $N \in \NN$.  Since $\tau$ was arbitrary it follows that we can cover $[0,\omega(\overline\zeta))$ with finitely many of these intervals of length $2 \eps / 3$, which proves that the flow is topologically mixing.
	\end{proof}
	
	\begin{remark}
		The proofs rely on being able to shadow a point via a synchronizing word and having a degree of control on Birkhoff sums from the Walters property.  Nothing in the proofs requires that the alphabet $A$ be finite nor does it require compactness.  The proof holds in the setting of transitive shifts with countable alphabets that possess a synchronizing word with a roof function satisfying the Walters property. Note that a Walters roof function on a countable state shift can be unbounded. In particular,  locally H\"older functions from \cite{Sarig99} can be unbounded, but they satisfy the Walters property.  
	\end{remark}

	\begin{remark}
		The proof of Theorem \ref{thm:no-mix-dense} does require compactness, so we cannot say anything about countable alphabets in that case.
	\end{remark}

	\section{Suspensions flows over subshifts of finite type}\label{sec:SFT}	
	
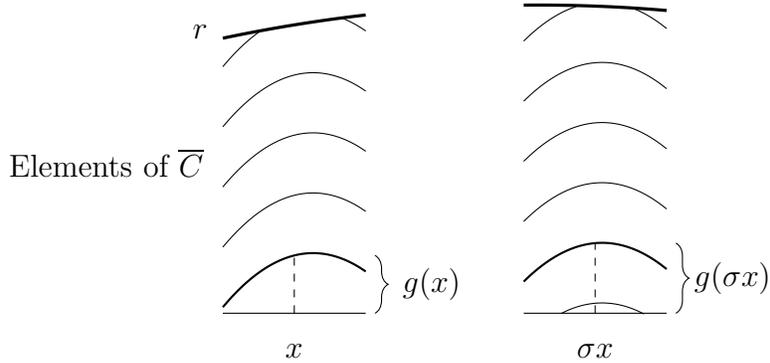
\begin{figure}[h!]
		\begin{tikzpicture}
			
			\def\dlta{.8}
			\def\scl{.5}
			
			\node at ({-1.75 + 2*\scl},2) {Elements of $\overline{C}$};
			\node at (.5,3.75) {$r$};
			
			\draw[scale=1, domain={1.6*\scl}:{5.4*\scl}, smooth, very thick, variable=\x] plot ({\x}, {-.025*(\x - 5)*(\x - 5) + 4.1});
			\draw[scale=1, domain={1.6*\scl+4}:{5.4*\scl+4}, smooth, very thick, variable=\x] plot ({\x}, {-.025*(\x - 5)*(\x - 5) + 4.1});
			
			\draw ({1.6*\scl},0) -- ({5.4*\scl},0);
			\draw ({1.6*\scl+4},0) -- ({5.4*\scl+4},0);
			\node at ({3.5*\scl},-.5) {$x$};
			\node at ({3.5*\scl+4},-.5) {$\sigma x$};
			
			\draw[scale=1, domain={1.6*\scl}:{5.4*\scl}, smooth, thick, variable=\x, name path=A] plot ({\x}, {-.5*(\x - 4*\scl)*(\x - 4*\scl) + 1*\dlta});
			\draw[scale=1, domain={1.6*\scl}:{5.4*\scl}, smooth, variable=\x] plot ({\x}, {-.5*(\x - 4*\scl)*(\x - 4*\scl) + 2*\dlta});
			\draw[scale=1, domain={1.6*\scl}:{5.4*\scl}, smooth, variable=\x] plot ({\x}, {-.5*(\x - 4*\scl)*(\x - 4*\scl) + 3*\dlta});
			\draw[scale=1, domain={1.6*\scl}:{5.4*\scl}, smooth, variable=\x] plot ({\x}, {-.5*(\x - 4*\scl)*(\x - 4*\scl) + 4*\dlta});
			\draw[scale=1, domain={1.6*\scl}:{2.6*\scl}, smooth, variable=\x] plot ({\x}, {-.5*(\x - 4*\scl)*(\x - 4*\scl) + 5*\dlta});
			\draw[scale=1, domain={4.8*\scl}:{5.4*\scl}, smooth, variable=\x] plot ({\x}, {-.5*(\x - 4*\scl)*(\x - 4*\scl) + 5*\dlta});
			
			\draw[scale=1, domain={2.6*\scl+4}:{4.8*\scl+4}, smooth, variable=\x] plot ({\x}, {-.5*(\x - 4*\scl-4)*(\x - 4*\scl-4) + 5*\dlta +.025*(\x - 9)*(\x - 9) - 4.1});
			\draw[scale=1, domain={1.6*\scl+4}:{5.4*\scl+4}, smooth, thick, variable=\x] plot ({\x}, {-.5*(\x - 4*\scl-4)*(\x - 4*\scl-4) + 5*\dlta +.025*(\x - 9)*(\x - 9) - 4.1 + \dlta});
			\draw[scale=1, domain={1.6*\scl+4}:{5.4*\scl+4}, smooth, variable=\x] plot ({\x}, {-.5*(\x - 4*\scl-4)*(\x - 4*\scl-4) + 5*\dlta +.025*(\x - 9)*(\x - 9) - 4.1 + 2*\dlta});
			\draw[scale=1, domain={1.6*\scl+4}:{5.4*\scl+4}, smooth, variable=\x] plot ({\x}, {-.5*(\x - 4*\scl-4)*(\x - 4*\scl-4) + 5*\dlta +.025*(\x - 9)*(\x - 9) - 4.1 + 3*\dlta});
			\draw[scale=1, domain={1.6*\scl+4}:{5.4*\scl+4}, smooth, variable=\x] plot ({\x}, {-.5*(\x - 4*\scl-4)*(\x - 4*\scl-4) + 5*\dlta +.025*(\x - 9)*(\x - 9) - 4.1 + 4*\dlta});
			\draw[scale=1, domain={1.6*\scl+4}:{3*\scl+4}, smooth, variable=\x] plot ({\x}, {-.5*(\x - 4*\scl-4)*(\x - 4*\scl-4) + 5*\dlta +.025*(\x - 9)*(\x - 9) - 4.1 + 5*\dlta});
			\draw[scale=1, domain={4.4*\scl+4}:{5.4*\scl+4}, smooth, variable=\x] plot ({\x}, {-.5*(\x - 4*\scl-4)*(\x - 4*\scl-4) + 5*\dlta +.025*(\x - 9)*(\x - 9) - 4.1 + 5*\dlta});
			
			\draw[dashed] ({3.5*\scl},0) -- ({3.5*\scl}, {-.5*(3.5*\scl - 4*\scl)*(3.5*\scl - 4*\scl) + 1*\dlta});			
			\draw[dashed] ({3.5*\scl+4},0) -- ({3.5*\scl+4}, {-.5*(3.5*\scl+4 - 4*\scl-4)*(3.5*\scl+4 - 4*\scl-4) + 5*\dlta +.025*(3.5*\scl+4 - 9)*(3.5*\scl+4 - 9) - 4.1 +\dlta});
			
			\draw [decorate,decoration={brace,amplitude=5pt,mirror,raise=0.4ex}]
			({5.5*\scl},0) -- ({5.5*\scl}, {-.5*(3.5*\scl - 4*\scl)*(3.5*\scl - 4*\scl) + 1*\dlta}) node[midway,xshift=2em]{$g(x)$};
			\draw [decorate,decoration={brace,amplitude=5pt,mirror,raise=0.4ex}]
			({5.5*\scl + 4},0) -- ({5.5*\scl+4}, {-.5*(3.5*\scl+4 - 4*\scl-4)*(3.5*\scl+4 - 4*\scl-4) + 5*\dlta +.025*(3.5*\scl+4 - 9)*(3.5*\scl+4 - 9) - 4.1 +\dlta}) node[midway,xshift=2em]{$g(\sigma x)$};
		\end{tikzpicture}
		\caption{Visualizing the transfer function in the phase space of the flow.  The thick curves of $\overline{C}$ give the value of $g$.}\label{fig:transfer-function}
	\end{figure}

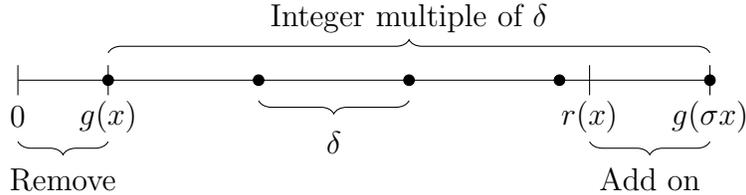
\begin{figure}[h!]
		\def\hash{0.2}
		\begin{tikzpicture}
			\draw (0,0) -- (9.2,0);
			\draw (0,{-\hash}) -- (0,{\hash});
			\draw (1.2,{-\hash}) -- (1.2,{\hash});
			\draw (7.6,{-\hash}) -- (7.6,{\hash});
			\draw (9.2,{-\hash}) -- (9.2,{\hash});
			
			\filldraw (1.2,0) circle (2pt);
			\filldraw (1.2+2,0) circle (2pt);
			\filldraw (1.2+4,0) circle (2pt);
			\filldraw (1.2+6,0) circle (2pt);
			\filldraw (1.2+8,0) circle (2pt);
			
			\node at (0,{-2.4*\hash}) {$0$};
			\node at (1.2,{-2.4*\hash}) {$g(x)$};
			\node at (7.6,{-2.4*\hash}) {$r(x)$};
			\node at (9.2,{-2.4*\hash}) {$g(\sigma x)$};
			
			\draw [decorate,decoration={brace,amplitude=5pt,mirror,raise=4.5ex}]
			(0,0) -- (1.2,0) node[midway,yshift=-3.2em]{Remove};
			\draw [decorate,decoration={brace,amplitude=5pt,raise=2ex}]
			(1.2,0) -- (9.2,0) node[midway,yshift=2em]{Integer multiple of $\delta$};
			\draw [decorate,decoration={brace,amplitude=5pt,mirror,raise=4.5ex}]
			(7.6,0) -- (9.2,0) node[midway,yshift=-3.2em]{Add on};
			\draw [decorate,decoration={brace,amplitude=5pt,mirror,raise=1.5ex}]
			(3.2,0) -- (5.2,0) node[midway,yshift=-2em]{$\delta$};
		\end{tikzpicture}
		\caption{Visualizing $r(x) - g(x) + g(\sigma x)$ along the orbit of $(x,0)$.  The dots represent elements of $\overline{C}$ in the orbit of $(x,0)$.}\label{fig:orbit-decomp}
	\end{figure}

	Here we present the proof of Theorem \ref{thm:SFT-cohomologous}.
	\begin{proof}
		Since the suspension flow $\ph$ is not topologically mixing, there exists a $0 < \delta < \min r/3$ satisfying the condition in equation \eqref{eqn:no-mix-orbit-length} for all periodic points $x \in X$ by Theorem \ref{thm:mix-dich} and Corollary \ref{cor:SFT-commensurate}.  
		
		Since $X$ is transitive, there exists a point $x_0$ that has a dense orbit in $X$.  Consider the collection of points in the phase space, $M$ of the flow
		\[
		C = \bigcup_{t \in \ZZ} \ph^{t\delta}(x_0,0).
		\] 
		We will show that $\overline{C}$ is the union of continuous ``curves" in $M$.
		
		Let $\eps > 0$ be small, and choose $k$ large enough to use the Walters property for this choice of $\eps$.  There exists an $\ell \in \NN$ such that for any $x\in X$ with $d(x,\sigma^n x) \leq 2^{-\ell}$, there exists a periodic point $p = \sigma^np$ such that $d(\sigma^j p, \sigma^j x) \leq 2^{-k}$ for all $0 \leq j < n$ by the Closing Lemma (see Theorem \ref{lem:closing}).   
		
		Let $y$ be any point in the orbit of $x_0$, and suppose that $d(\sigma^n y,y)\leq 2^{-\ell}$ for some $n \in \NN$.  Let $z \in X$ be the periodic point satisfying $d(\sigma^j z, \sigma^jy) \leq 2^{-k}$ for all $0 \leq j < n$.  
		
		Thus, we have
		\begin{equation}
			\left|\sum_{j=0}^{n-1} r(\sigma^j y) - r(\sigma^j z) \right| < \eps
		\end{equation}
		by the Walters property.
		
		Since $z$ has period $n$, there is a $K \in \NN$ such that 
		\begin{equation}\label{eqn:closing+Wal}
			\left|\sum_{j=0}^{n-1} r(\sigma^j y) - K\delta \right| < \eps,
		\end{equation}
		by Proposition \ref{prop:mult-sync}.
		
		If we let
		\[
		\tau = \sum_{j=0}^{n-1} r(\sigma^j y),
		\]
		then $\ph^{\tau}(y,0) = (\sigma^n y,0)$. Moreover,
		\begin{equation}\label{eqn:closing+Walters-a}
			\ph^{K\delta}(y,0) \in \ph^{(\tau-\eps,\tau+\eps)}(y,0) = \ph^{(-\eps,\eps)}(\sigma^ny,0),
		\end{equation}
		by equation \eqref{eqn:closing+Wal}.  If $(y,a) \in C$ for some $\eps < a < r(\sigma^n y) - \eps$, then $\ph^{K\delta}(y,a) \in C$ by construction.  Equation \eqref{eqn:closing+Walters-a} implies that
		\[
		\ph^{K\delta}(y,a) \in  \ph^{(\tau-\eps,\tau+\eps)}(y,a) = \ph^{(a-\eps,a+\eps)}(\sigma^ny,0) = \ph^{(-\eps,\eps)}(\sigma^ny,a),
		\]
		so the real coordinates of $\ph^{K\delta}(y,a)$ and $\ph^{\tau}(y,a)$ are within $\eps$ of each other in the flow direction.  Furthermore, because $a \in(\eps , r(\sigma^ny) - \eps)$ the real coordinate of $\ph^{K\delta}(y,a)$ must be within $\eps$ of $a$, which is the starting height of $(y,a) \in C$.  
		
		Additionally, if $(\sigma^n y,b) \in C$ and $\eps < b < r(y) - \eps$, then there must be some $j \in \ZZ$ and a point $(y,a + j\delta) \in C$ such that $\ph^{K\delta}(y,a+j\delta) = (\sigma^n y,b)$ and $|b - (a+j\delta)| < \eps$. 
		
		Since this is true for any return to the $2^{-\ell}$-neighborhood of $y$ and $\eps$ was arbitrary, it follows that $C$ is a collection of continuous curves on the orbit of $x_0$.  Since the orbit of $x_0$ is dense, $\overline{C}$ is a collection of continuous curves in $M$, and $\ph^{\delta}(\overline{C}) = \overline{C}$.  
		
		Given a word $w \in \LLL$, the intersection $([w] \times \RR) \cap \overline{C}$ is a collection of continuous curves that are all separated by a distance of $\delta$ is the flow direction.  We say that a curve $\gamma$ in $([w] \times \RR) \cap \overline{C}$ sits over all of $[w]$ if $\gamma \cap (\{x\}\times\RR) \neq \emptyset$ for all $x \in [w]$.  It is possible for a curve of $([w] \times \RR) \cap \overline{C}$ to not sit over all of $[w]$ if it is too close to the base or the roof of the suspension.
		
		Since $\delta < \min r/3$, we can fix $N$ sufficiently large so that for all $w \in \LLL_N$, we have that $([w] \times \RR) \cap \overline{C}$ contains at least one continuous curve that sits over all of $[w]$, as sketched in Figure \ref{fig:transfer-function}.  For all $w \in \LLL_N$, let $\gamma_w$ be the curve closest to the base that sits over all of $[w]$.
		
		Define $g\colon X \to (0,\infty)$ by setting $g(x)$ to be the height of the curve $\gamma_w$ over $x \in [w]$.  The function $g$ is well-defined and continuous because, for all $w \in \LLL_N$, we chose each $\gamma_w$ to be a continuous curve that sits over all of $[w]$.  
		
		Now, consider the function
		\begin{equation}
			s(x) = r(x) - g(x) + g(\sigma x).
		\end{equation}
		We claim that $s(x)$ is locally constant. Removing $g(x)$ from $r(x)$ moves an element of $\overline{C} \cap (\{x\} \times \RR)$ to the base.  By adding $g(\sigma x)$, we force $r(x) - g(x) + g(\sigma x)$ to have a value in $\delta \ZZ$ because the time it takes to flow from $(x,g(x))$ to $(\sigma x, g(\sigma x))$ must be an integer multiple of $\delta$, as $(x,g(x))$ and $(\sigma x, g(\sigma x))$ are elements of $\overline{C}$, see Figure \ref{fig:orbit-decomp}.  
		
		Since $s$ is continuous and maps into $\delta \ZZ$, $s$ must be locally constant.  Therefore, $r$ is cohomologous to a locally constant roof function that only takes values in $\delta \ZZ$.
	\end{proof}
	
	One can easily verify the following lemma as a consequence of Lemma \ref{lem:fin-range}.
	\begin{lemma}\label{lem:countable-roofs}
		If $X$ is a shift over a finite alphabet, then there are countably many locally constant roof functions $r \colon X \to \NN$.
	\end{lemma}	
	Let $W(X)$ denote the set of roof functions of $X$ that satisfy the Walters property.  Given $r \in W(X)$ define
	\begin{equation}\label{eqn:Walters-norm}
		\| r \|_{\text{Wal}} = 2\| r\|_{\infty} + \sup_{m \geq 1} \max_{d(x,y) < 2^{-m}} \left| \sum_{j=0}^{m-1} r(\sigma^j x) - r(\sigma^j y) \right|.
	\end{equation}
	As the notation suggests, $\|\cdot\|_{\text{Wal}}$ is a norm, and $W(X)$ is a Banach space when it is equipped with this norm  \cite{tB01}. 
	
	We now prove Theorem \ref{thm:SFT-G-delta}.
	\begin{proof}
		Let 
		\[W^* = \{r \in W(X) : r \text{ does not produce a topologically mixing flow over } X\}.\]
		Suppose $r \in W^*$, and let $W_r = \{s \in W(X) : s \text{ is cohomologous to } r \}$.  $W_r$ is a closed subset of $W(X)$ because $W(X)$ is a Banach space and by the Liv\v sic theorem (see Theorem \ref{thm:Livsic}).
		
		Now let $W^*_r = \{ks : k \in (0,\infty) \text{ and } s \in W_r\}$, which is also closed.  There exists a $\delta > 0$ and an $s \in W_r$ such that $s \colon X \to \delta \NN$ by Theorem \ref{thm:SFT-cohomologous}. This implies that $s/\delta \in W^*_r$.
		
		Since $W^*_r$ contains a roof function that maps into $\NN$ and there are only countably many of these roof functions by Lemma \ref{lem:countable-roofs},  $W^*(X)$ can be written as a countable union of closed sets.  Clearly, $(W^*)^C$ is dense in $W(X)$ and can be expressed as a countable intersection of open sets.  Therefore, $(W^*)^C$ is a dense, $G_\delta$ subset of $W(X)$.
	\end{proof}

	\section{Suspension flows over \texorpdfstring{$\beta$}{beta}-shifts} \label{sec:beta-shifts}
	
	$\beta$-shifts arise naturally from codings of interval expanding maps $T_\beta \colon [0,1]\to[0,1]$ defined by $T_\beta(x) = \beta x \mod 1$, and they have connections to $\beta$-expansions of real numbers.  For $x \in \RR$, we let $\lfloor x \rfloor$ denote the integer part of $x$.  If $\beta > 1$, then the $\beta$-expansion of a number $x \in [0,1]$ is a sequence $\{a_n\}_{n=1}^\infty$ such that $a_n \in \{0,1,\dots, \lfloor \beta \rfloor\}$ for all $n \in \NN$ and
	\[
	a_n = \lfloor \beta T_{\beta}^{n-1}(x) \rfloor.
	\]
	\begin{definition}
		The closure of the set of all $\beta$-expansions of $x \in [0,1]$ is a \emph{$\beta$-shift}.
	\end{definition}
	
	If $X \subset \{0,\dots,b\}^\NN$ is a $\beta$-shift, then there is a point $\nu(\beta) \in X$ such that for any $x \in X$, $\sigma^j x \preceq \nu(\beta)$ for all $j\geq 0$ where $\preceq$ is denotes the lexicographical ordering on $\{0,\dots,b\}^\NN$.  $\nu(\beta)$ is the $\beta$-expansion of 1 \cite{wP60}.  Although $\beta$-shifts are typically studied as a one-sided shift, one can extend a one-sided shift to a two-sided shift in the usual way.  In the context of a suspension flow, it is more natural for the dynamics in the base to be a homeomorphism, so we will work with the two-sided version of the $\beta$-shift.  
	
	A one-sided $\beta$-shift can be understood by comparing it to walks on a countable state directed graph with multi-edges.  An edge-walk on the graph corresponds to points in the $\beta$-shift.  We label the vertices $V_1,V_2,\dots$ and designate $V_1$ as the \emph{initial vertex}. There are two types of edges in the graph.  There is a single edge that connects vertex $V_n$ to vertex $V_{n+1}$ for all $n\in \NN$.  There can also be multiple edges from $V_n$ to the initial vertex $V_1$; however, these edges may only be present for some $n\in \NN$.
	
\begin{figure}[h!]
		\centering
		\def\hgap{2.3}
		\def\lilspace{0.75}
		\begin{tikzpicture}
			\tikzset{vertex/.style = {shape=circle,draw,minimum size=1.5em}}
			\tikzset{edge/.style = {->,> = latex'}}
			
			\node[vertex] (1) at (0,0) {$V_1$};
			\node[vertex] (2) at ({\hgap},0) {$V_2$};
			\node[vertex] (3) at ({2*\hgap},0) {$V_3$};
			\node[vertex] (4) at ({3*\hgap},0) {$V_4$};
			\node[vertex] (5) at ({4*\hgap},0) {$V_5$};
			
			\draw[edge, shorten >={\lilspace}] (1) to node[above]{2} (2);
			\draw[edge, shorten >={\lilspace}] (2) to node[above]{0} (3);
			\draw[edge, shorten >={\lilspace}] (3) to node[above]{1} (4);
			\draw[edge, shorten >={\lilspace}] (4) to node[above]{2} (5);
			
			\draw[edge, shorten >={\lilspace}] (1.130) .. controls (-1,.5) and (-1,-.5) ..  node[left]{0} (1.250);
			\draw[edge, shorten >={\lilspace}] (1.120) .. controls (-1.75,1) and (-1.75,-1) ..  node[left]{1} (1.270);
			
			\draw[edge, shorten >={\lilspace}, bend left = 15] (3) to node[below]{0} (1.290);
			\draw[edge, shorten >={\lilspace}, bend left = 27] (4) to node[below]{1} (1.280);	
			\draw[edge, shorten >={\lilspace}, bend left = 45] (4) to node[below]{0} (1.270);	
			
			\draw[edge,dashed] (5) to ({5*\hgap-.5},0);
			
		\end{tikzpicture}
		\caption{Representation of $\beta$ shift as a countable state shift}\label{fig:beta-graph}
	\end{figure}
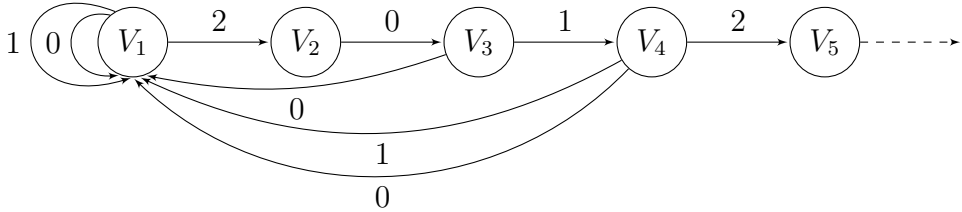

	The label for the edge $V_n$ to $V_{n+1}$ is $\nu(\beta)_n$.  If the label from $V_n$ to $V_{n+1}$ is greater than 0, then there are edges from $V_n$ to the initial vertex $V_1$ for every label less than $\nu(\beta)_n$.  That is, if $\nu(\beta)_n \geq 1$, then there are $\nu(\beta)_n-1$ edges from $V_n$ to $V_1$ with labels $\{0,\dots, \nu(\beta)_n-1\}$.  For example, in Figure \ref{fig:beta-graph}, the label from $V_4$ to $V_5$ is $2$, so there are edges label $0$ and $1$ from $V_4$ to $V_1$.  To obtain the two-sided $\beta$-shift, we consider bi-infinite walks on this graph.
	
	A $\beta$-shift is sofic if and only if $\nu(\beta)$ is eventually periodic.  In this case, the $\beta$-shift is an irreducible sofic shift because it can be presented by a strongly connected directed graph on finitely many vertices \cite{fischer}.  
	
	A $\beta$-shift has specification if and only if $\nu(\beta)$ does not possess arbitrarily long blocks of $0$s \cite{CT12}. Both irreducible sofic shifts and shifts with specification possess synchronizing words \cite{bertrand,jonoska}, so Theorem \ref{thm:mix-dich} applies to these shifts. However, there are $\beta$-shifts that do not possess synchronizing words; these shifts possess points with arbitrarily long strings of 0s.  Additionally, the $\beta \in \RR$ that produce sofic shifts or shifts with specification is a set with zero Lebesgue measure \cite{schmeling}.  We prove Theorem \ref{thm:mix-dich-beta} in the case when $\nu(\beta)$ can have arbitrarily long strings of 0s, so the shift does not possess a synchronizing word.  More information about $\beta$-shifts can be found in \cite{CT12,wP60,Renyi,kT05}
	
	The proofs of Theorem \ref{thm:mix-dich-beta} for $\beta$-shifts are similar to those of Theorem \ref{thm:mix-dich}. 
	\begin{proof}
		
		We begin by proving the necessary condition of Theorem \ref{thm:mix-dich-beta}.  Fix $0 < \eps < \delta/16$, and let $k \in \NN$ be the constant needed to leverage the Walters property for this choice of $\eps$.  Consider the cylinder $[0^{2k+1}]$, and suppose $y \in [0^{2k+1}]$ such that $\sigma^n y \in [0^{2k+1}]$ for any $n \gg k$.  Then 
		\[y = \dots 0^k . 0^{k+1} w 0^{2k+1} \dots \]
		for some $w \in \LLL_{n-2k-1}$.  The word $0^{2k+1}w0^{2k+1}$ corresponds to some walk on the countable state graph.  Thus, there is some word $v$ such that $v0^{2k+1}w0^{2k+1}$ codes a walk that starts at the initial vertex.  Hence, there is a periodic point $p = \overline{0^{k}.0^{k+1}w0^{2k+1}0^{\ell}v}$, where the ``." indicates the center position of $p$ so that $p \in [0^{2k+1}]$.  We choose $\ell$ sufficiently large so that  $0^{2k+1} w 0^{2k+1} 0^{\ell}$  corresponds to a walk that returns to the initial vertex and $v 0^{2k+1} 0^{\ell}$ contains a subword that corresponds to a walk that starts and ends at the initial vertex $V_1$.
		
		That is, let $y^*$ be the one-sided tail of $y$ so that we may write 
		\[
		y = \dots 0^k.0^{k+1} w 0^{2k+1} y^*.
		\]  If $y^*$ contains a nonzero symbol, then there is an $\ell'$ so that $y^* = 0^{\ell'-1}\dots$ where $\ell'$ is the maximal number of $0$s before $y^*$ has a nonzero symbol.  Thus for $\ell \geq \ell'$, $0^{2k+1} w 0^{2k+1+\ell}$ corresponds to a return to the initial vertex (or contains such a subword).  If $y^*$ is all 0 symbols, then $0^{2k+1}w0^{2k+1}$ already corresponds to a return to the initial vertex (or contains such a subword), and we may choose $\ell \geq 0$.
		
		Likewise, $v0^{2k+1}$ is a walk that begins at the initial vertex by the choice of $v$.  There must be an $\ell''$ such that $v0^{2k+1+\ell''}$ is a walk that begins and ends at the initial vertex.  We may choose $\ell > \ell', \ell''$, as we may add on arbitrarily many 0s since the walk is at the initial vertex after $\ell'$ or $\ell''$ steps.  This implies that there exists a periodic point $q = \overline{0^k.0^{k+1+\ell}v}$.  
		
		Note that $n = 2k + 1 + |w|$ we have that
		\begin{equation}\label{eqn:p-delta}
			\omega(p) = \sum_{j=0}^{n+2k+\ell+|v|} r(\sigma^j p) = C_1\delta
		\end{equation}
		for some $C_1 \in \NN$ because $p$ is a periodic point in $[0]$.  Moreover, $p_{[-k,n+k]} = y_{[-k,n+k]}$, so
		\begin{equation}\label{eqn:shadow-p-y}
			\left|\sum_{j=0}^{n-1} r(\sigma^j y) - r(\sigma^j p) \right| < \eps
		\end{equation}
		by the Walters property.
		
		Likewise, since $q \in [0]$ and is periodic, it follows from the hypothesis that there is a $C_2 \in \NN$ such that
		\begin{equation}\label{eqn:q-delta}
			\omega(q) = \sum_{j=0}^{2k+\ell +|v|} r(\sigma^j q) = C_2\delta.
		\end{equation}
		Additionally, observe that $q_{[-k,3k+l+|v|+2]} = \sigma^{n} p_{[-k,3k+l+|v|+2]}$, so
		\begin{equation}\label{eqn:shadow-q-p}
			\left|\sum_{j=0}^{2k+\ell+|v|} r(\sigma^j q) - r(\sigma^{j+n} p)\right| < \eps
		\end{equation}
		by the Walters property.
		
		Using equations \eqref{eqn:p-delta} and \eqref{eqn:q-delta} followed by the triangle inequality and \eqref{eqn:shadow-p-y} and  \eqref{eqn:shadow-q-p}, we get
		\begin{align*}
			\left|(C_1-C_2)\delta - \sum_{j=0}^{n-1} r(\sigma^jy) \right| & =  \left|\omega(p) - \omega(q) - \sum_{j=0}^{n-1} r(\sigma^jy) \right| \\ &\leq \left| \sum_{j=0}^{n-1} r(\sigma^j p) - r(\sigma^j y) \right| + \left|\sum_{j=0}^{2k+\ell+|v|-1} r(\sigma^j q) - r(\sigma^{j+n} p)\right|\\
			& < 2\eps
		\end{align*}
		By the same reasoning as in the proof of Theorem \ref{thm:mix-dich}, the flow is not topologically mixing.

		We now prove the sufficient condition by way of contraposition.  The proof is very similar to that of Theorem \ref{thm:mix-dich} in \S\ref{sec:mix-dict-pf}.

		Let $u,w \in \LLL$ be any words with $|u| = |w|$.  Let $U = [u] \times (0,\eps)$ and $W = [w]\times (0,\eps)$.  Fix $0< \gamma < \eps/21$.  There exists a $k$ such that if $x,y\in X$ satisfying $x_{[-k,n+k]}=y_{[-k,n+k]}$, then $|\sum_{j=0}^{n-1}r(\sigma^j x) - \sum_{j=0}^{n-1}r(\sigma^j y)|<\gamma$ by the Walters property.
		
		As before, for any periodic point $p \in X$, let $\omega(p)$ denote the period of $(p,0)$ under the flow.  By hypothesis, there exists periodic points $p,q \in [0]$ such that either $\omega(p)$ and $\omega(q)$ are incommensurate or $\gcd(\omega(p),\omega(q)) \leq \gamma$.  Consequently, there exists $\alpha,\beta \in \NN$ such that 
		\[0 < |\alpha\omega(p)-\beta\omega(q)| \leq \gamma.\]
		
		Since $p,q \in [0]$ are periodic points there exists word $v_1,v_2 \in \LLL$ such that $p = \overline{v_1} = \overline{0v'}$ and $q = \overline{v_2} = \overline{0v''}$.  Note that here we are slightly abusing notation as $v_1$ and $v_2$ are words and do not represent symbols of some word $v$.  Since $p$ and $q$ are periodic points and $\nu(\beta)$ is neither periodic nor eventually periodic, we know that $0v'$ and $0v''$ are words that correspond to cycles on the countable state shift associated to $X$.  This means that the walks associated to these periodic points must visit the initial vertex $V_1$ of the graph since all cycles of the graph must go through $V_1$.
		
		Without loss of generality, we may assume that the first symbol in $v_1$ and $v_2$ is $0$ and corresponds to an edge in the graph that ends at the initial vertex.  Hence, $v_1^nv_2^m \in \LLL$.
		
		Additionally, there is a word $w_f$ such that $w w_f \in \LLL$ and $w_f$ corresponds to a position on the graph that can be followed by a ``fall" to the initial vertex.  In other words, there is a symbol $a$ such that $a>0$ and $ww_fa \in \LLL$.  Similarly, there are words $u_c$ and $u_f$ such that $u_c u u_f \in \LLL$ and $u_c$ is the climb from the initial vertex to $u$ and $u_f$ is the word after $u$ that corresponds to a return back to the initial vertex.  Unlike $w_f$, we will allow the last symbol of $u_f$ to correspond to the edge back to the initial vertex.  Note that it is possible that $w_f$, $u_c$, or $u_f$ is the empty word.  Now define $\zeta = u_c u u_f$, and observe that $\overline{\zeta} \in X$.
		
		Let
		\[
		x_{n,m} = \dots w^* w w_f v_1^n v_2^m  \zeta^\infty \in [w]
		\]
		for some permissible $w^*\in\LLL_{\geq k}$.
		
		Fix $n_0 = m_0 = 3k$ and consider $x^* \coloneqq x_{n_0,m_0}$.  Now let $n,m \geq 3k$.  By shadowing the orbit of $x_{n,m}$ with the appropriate points, we can repeatedly apply the Walters property just like in the proof of Theorem \ref{thm:mix-dich}.
		\begin{equation}\label{eqn:beta-first-est}
			\left|\sum_{j=0}^{|ww_f|+k|v_1|-1} r(\sigma^j x_{n,m}) - \sum_{j=0}^{|ww_f|+k|v_1|-1} r(\sigma^j x^*) \right| < \gamma
		\end{equation}
		Let $x_{n,m}^{(1)} = \sigma^{|ww_f|+k|v_1|}x_{n,m}$.
		\begin{equation}\label{eqn:beta-shadow-p}
			\left|\sum_{j=0}^{(n-2k)|v_1|-1} r(\sigma^{j} x_{n,m}^{(1)}) - \sum_{j=0}^{(n-2k)|v_1|-1} r(\sigma^j p) \right| < \gamma
		\end{equation}
		Let $x_{n,m}^{(2)} = \sigma^{(n-2k)|v_1|}x_{n,m}^{(1)}$.
		\begin{equation}\label{eqn:beta-shadow-second}
			\left|\sum_{j=0}^{k(|v_1|+|v_2|)-1} r(\sigma^{j} x_{n,m}^{(2)}) - \sum_{j=0}^{k(|v_1|+|v_2|)-1} r(\sigma^{j+|ww_f|+2k|v_1|} x^*) \right| < \gamma
		\end{equation}
		Let $x_{n,m}^{(3)} = \sigma^{k(|v_1|+|v_2|)}x_{n,m}^{(2)}$.
		\begin{equation}\label{eqn:beta-shadow-q}
			\left|\sum_{j=0}^{(m-2k)|v_2|-1} r(\sigma^{j} x_{n,m}^{(3)}) - \sum_{j=0}^{(m-2k)|v_2|-1} r(\sigma^{j} q) \right| < \gamma
		\end{equation}
		Let $x_{n,m}^{(4)} = \sigma^{(m-2k)|v_2|}x_{n,m}^{(3)}$.
		\begin{equation}\label{eqn:beta-arrive-est}
			\left|\sum_{j=0}^{k(|v_2|+|\zeta|)-1} r(\sigma^{j} x_{n,m}^{(4)}) - \sum_{j=0}^{k(|v_2|+|\zeta|)-1} r(\sigma^{j+|ww_f|+3k|v_1|+2k|v_2|} x^*) \right| < \gamma
		\end{equation}
		Let $x_{n,m}^{(5)} = \sigma^{k(|v_2|+|\zeta|)}x_{n,m}^{(4)}$ and $N\in\NN$.
		\begin{equation}\label{eqn:beta-return-to-zeta-time}
			\left|\sum_{j=0}^{N-1} r(\sigma^{j} x_{n,m}^{(5)}) - \sum_{j=0}^{N-1} r(\sigma^{j} \overline{\zeta}) \right| < \gamma.
		\end{equation}
		Recall, $\omega(\overline{\zeta}) = \sum_{j=0}^{|\zeta|-1}r(\sigma^j\overline{\zeta})$.  By equation \eqref{eqn:beta-return-to-zeta-time}, for any $n,m\geq 3k$ we have
		\[\left| \sum_{j=0}^{N|\zeta|-1}r(\sigma^jx_{n,m}^{(5)}) -N\omega(\overline{\zeta}) \right| < \gamma\]
		This implies that the return time of $x_{n,m}^{(5)}$ to $U$ is in  $(\omega(\overline{\zeta}) - \gamma, \omega(\overline{\zeta}) + \gamma) $. 
		
		Let 
		\begin{align*}
			K & = \sum_{j=0}^{|ww_f|+k|v_1|-1} r(\sigma^j x^*) + \sum_{j=0}^{k(|v_1|+|v_2|)-1} r(\sigma^{j+|ww_f|+2k|v_1|} x^*)\\
			& \qquad + \sum_{j=0}^{k(|v_2|+|\zeta|)-1} r(\sigma^{j+|ww_f|+3k|v_1|+2k|v_2|} x^*)
		\end{align*}
		
		Let $T_{n,m} = |ww_f v_1^nv_2^m \zeta^k|$.  If we use equations \eqref{eqn:beta-first-est}-\eqref{eqn:beta-arrive-est} and apply the triangle inequality, we have
		\begin{equation}\label{eqn:beta-full-arrive-est}
			\left| \sum_{j=0}^{T_{n,m}-1}r(\sigma^jx_{n,m}) - (K + (n-2k)\omega(p) + (m-2k)\omega(q))\right| < 5\gamma
		\end{equation}
		for any pair $n,m$.  
		
		Let $\tau \in (0,\omega(\overline{\zeta}))$ be arbitrary.  By Lemma \ref{lem:gcd-dense} we know that there exists $n,m, C \in \NN$ such that 
		
		\[
		\left| \tau +C\omega(\overline{\zeta}) - (K + (n-2k)\omega(p) + (m-2k)\omega(q)) \right|  < \gamma.
		\]
		Applying equation \eqref{eqn:beta-full-arrive-est} we have
		\[\left| \tau + C\omega(\overline{\zeta}) - \sum_{j=0}^{T_{n,m}-1}r(\sigma^jx_{n,m}) \right| <  6\gamma\]
		and \eqref{eqn:beta-return-to-zeta-time} estimates the subsequent returns to $U$ by
		\[\left| \tau (C+N)\omega(\overline{\zeta}) - \sum_{j=0}^{T_{n,m}+N|\zeta|-1}r(\sigma^jx_{n,m}) \right|  < 7\gamma.\]
		
		For $N\in\NN$, the time $\sum_{j=0}^{T_{n,m}+N|\zeta|-1}r(\sigma^jx_{n,m})$ corresponds to the moment that the image of $x_{n,m}\times (0,\eps) \subset W$ under the flow is contained in $U$.  Because $\gamma < \eps/21$ we are guaranteed that $x_{n,m}\times(0,\eps)$ is passing through $U$ from $\tau + (C+N)\omega(\overline{\zeta}) - \eps/3$ to $\tau+(C+N)\omega(\overline{\zeta}) + \eps/3$ for all $N \in \NN$.
		
		In other words, if we let $T = \sum_{j=0}^{T_{n,m}-1}r(\sigma^jx_{n,m})$, then for all $t \in (\tau + (C+N)\omega(\overline\zeta) - \eps/3, \tau + (C+N)\omega(\overline\zeta) + \eps/3)$ we are guaranteed that $\ph^{T+t}(W)\cap U\neq\emptyset$ for all $N\in\NN$.  Since $\tau$ was arbitrary it follows that we can cover $[0,\omega(\overline\zeta))$ with finitely many of these $\eps$-intervals, which proves that the flow is topologically mixing.
	\end{proof}
	
	\section{Proof of Theorems \ref{thm:loc-const-com-nonmixing} and \ref{thm:no-mix-dense}}\label{sec:density}

	\begin{lemma}\label{lem:fin-range}
		If $X$ is a compact metric space and $r$ is a locally constant roof function, then the range of $r$ is finite.
	\end{lemma}
	\begin{proof}
		Let $x\in X$.  There exists an open neighborhood $U_x$ of $x$ such that $r$ is constant on $U_x$ by definition.  The collection of $\{U_x\}$ is an open cover of $X$.  Therefore, there is a finite subcover since $X$ is compact.  Since $r$ is constant on each subset of this subcover, we get that the range of $r$ must be finite.
	\end{proof}
	
	The following is the proof of Theorem \ref{thm:loc-const-com-nonmixing}
	\begin{proof}
		Let $\ph^t$ be the suspension flow over $X$ by the roof function $r$.  If $r$ is locally constant, then the range is finite by Lemma \ref{lem:fin-range}.  Let $\{r_1,\dots,r_n\}$ be the range of $r$ with $r_1 \leq \cdots \leq r_n$.  Since these values are all mutually commensurate, there exists a $\delta$ such that $r_i/\delta \in \NN$ for all $1\leq i \leq n$.  Without loss of generality, we can choose $\delta < r_1$.
		
		For $1\leq i\leq n$, let $X_i = r^{-1}(\{r_i\})$.  Note that $X_i$ is open because $r$ is locally constant, and $X_i$ is closed because it is the continuous preimage of a closed set.  Now define
		\[C = \bigcup_{i=1}^n \bigcup_{j=0}^{\frac{r_i}{\delta} - 1} X_i \times \{j\delta\}.\]
		We let $C$ be a transverse cross-section of the flow $\ph^t$ to obtain a new suspension flow which we will call $\psi^t$.  For all $1\leq i\leq n$ we have that $r_i$ is an integer multiple of $\delta$, so we get that the return time to $C$ under the flow is always $\delta$.  Thus $\psi^t$ is a suspension flow with a constant roof function of $\delta$.  Therefore, $\psi^t$ and $\ph^t$ are not topologically mixing.
	\end{proof}

	We now provide a formal proof of Theorem \ref{thm:no-mix-dense}.
	\begin{proof}
		Let $\epsilon>0$.  Since $X$ is compact and $r$ is continuous, it follows that $r$ is uniformly continuous.  Hence, there exists a $\delta>0$ such that for any $x,y\in X$ where $d(x,y)<\delta$, we get that $|r(x)-r(y)|<\epsilon$.   There exists an $n$ such that $2^{-n}<\delta$.  
		
		Suppose $w\in\LLL_{2n+1}$ and $x,y\in[w]$.  Since $|w|=2n+1$, it follows that $d(x,y)\leq 2^{-n-1}<\delta$, and $|r(x)-r(y)|<\epsilon$.  Since $\LLL_{2n+1}$ is a finite set, enumerate the words from 1 to $N$.  The collection of cylinders $\{[w_i]\}_{i=1}^N$ is a cover of $X$, and these cylinders are all mutually disjoint. 
		
		For each $w_i$ we may pick $x_i$ where $x_i\in[w_i]$ and let $k_i=r(x_i)$.  Define a roof function $s\colon X\to \RR$ by $s|_{[w_i]}=k_i$.  Clearly, $s$ is locally constant.
		
		For any $x\in X$, we have that $x\in[w_i]$ for some $w_i$.  This means that $d(x,x_i)<\delta$ and $|r(x)-k_i|=|r(x)-r(x_i)|<\epsilon$.  Since $x$ is arbitrary, this shows that the set of locally constant roof functions is dense in the set of all continuous roof functions. 
		
		Since the collection of locally constant roof functions with rational range is dense in the set of all locally constant roof functions, we are finished.
	\end{proof}

	\section{Topologically mixing shift with finitely many periodic points}\label{sec:mix-no-per}
	Here we present the construction of a topologically mixing shift space with only one fixed point and a period two orbit.  This construction is inspired by the structure of a weakly topologically mixing shift from a construction in \cite{FKKL} but echoes many of the ideas in \cite{Pet}.
	
	We will begin by constructing a one-sided subshift of $\{0,1\}^\NN$. However, we first require an element $a$ in $\{0,\dots,N\}^\NN$ with the following properties.
	\begin{enumerate}
		\item $a$ is almost periodic, but not periodic.
		\item $a_n\in\{3,4,\dots N\}$ for all $n\in\NN$
		\item For any $i,k\in \NN$ the word $a_{i}\dots a_{i+k}$ must appear an infinite number of times in $a$.\label{item:repeat}
	\end{enumerate}
	For reference, a point $x$ is an \textit{almost periodic point}, if for any neighborhood $U$ of $x$, there exists an $N \in  \NN$ such that $\{\sigma^{n+i}(x) : i=0,1,\dots, N\} \cap U \neq \emptyset$ for all $n\in\NN$.
	
	The first two items are important to ensure that the shift space we construct has the desired periodic structure.  The third item is true for the orbit of any almost periodic point by definition, but we highlight it here as it is needed to obtain topological mixing.  We know that such sequences exist as the Morse-Thue sequence satisfies all of these properties if we increase the terms in the sequence by 3.  Indeed, if $a$ were to satisfy all three of these conditions, then the closure of its orbit would be a minimal shift \cite{ML}.
	
	Let $u_1=101$, $u_2=10101$, $u_3=1010101$, and so on (once again, we are abusing notation here as $u_i$ represents a word, not a symbol of a word).  For $i\geq1$, define $E_i=\{u_1,u_2,\dots,u_i\} \cup \bigcup_{j\geq 2}1^j$ where $j \in \NN$.  For the almost periodic sequence $a$ satisfying the previously mentioned properties, define 
	\[B_i=\{x\in\{0,1\}^\NN : x=v_10^{a_1}v_20^{a_2}v_30^{a_3}\dots\text{ where }v_j\in E_i\}.\]
	Let
	\[J_i=\bigcup_{k=0}^\infty\sigma^k(B_i).\]
	Clearly $J_i\subseteq J_{i+1}$ for all $i\in\NN$.  Let $S=\overline{\bigcup_{i=1}^\infty J_i}$, which is clearly closed.  We must show that $\sigma(S) = S$.  By construction we have that $\sigma\left(\bigcup_{i=1}^\infty J_i\right) \subset \bigcup_{i=1}^\infty J_i$, so $\sigma(S)\subset S$ by continuity of $\sigma$.  We must show the opposite inclusion.
	
	For any $x \in \bigcup_{i=1}^\infty J_i$, we know that $x = \sigma^j y$ for some $j \in \NN$ and $y \in B_k$ for some $k$.  If $j \geq 1$, then we know that $x = \sigma(\sigma^{j-1}y) \in J_k$ and $\sigma^{j-1}y \in J_k$.  
	
	If $j=0$, then we know that $x \in B_k$ for some $k$ and begins with a word $v\in E_k$.  If $v = 10\dots101$, then we can find a point $y \in B_{k+1}$ such that $y = 10x$, so $x =\sigma^2y$.  Similarly, if $v=11$, then we can find $y\in B_k$ where $y = 1x$.  Hence $\sigma y = x$.  Therefore, $\bigcup_{i=1}^\infty J_i \subset \sigma(\bigcup_{i=1}^\infty J_i)$.
	
	Because $\bigcup_{i=1}^\infty J_i$ is dense in $S$, by continuity of $\sigma$ we obtain that $S \subset \sigma(S)$.  That is, for all $x\in S$, there exists $x'\in S$ such that $\sigma(x') = x$.  Hence $\sigma(S) = S$, so $S$ is invariant and $S$ is a subshift of $\{0,1\}^\NN$.  
	
	\begin{proposition}
		$S$ contains only two periodic orbits.
	\end{proposition}
	\begin{proof}
		Clearly, $1^\infty$ and $(01)^\infty$ are periodic points in $S$, since $\bigcup_{i=1}^\infty J_i$ contains points with arbitrarily long strings of 1s and 01 blocks.  We must show that $S$ contains no other periodic orbits.  
		
		Suppose $x \in S$ is a periodic point. Then there must be a word $w \in \LLL$ such that $x = w^\infty$.  Without loss of generality we may assume that $w = u_00^{a_i}u_1 \dots u_{k}0^{a_{i+k}}$ for some $i,k \in \NN$ and where $u_j$ are words of repeated $1$s or $10$s. 
		
		Since the blocks of $0$s in $w$ are always preceded and followed by a 1, the lengths of these blocks of $0$s cannot be altered by the $u_j$ words.   Because $w^n \in \LLL$ this would mean that the sequence $a$ has a block $(a_i\dots a_{i+k})^n$ contained in it.  Since we can find such a block for every $n \in \NN$, it would follow that $\overline{\bigcup_{j=1}^\infty \sigma^j a}$ has at least one periodic orbit.  However, the sequence $a$ was chosen to be an almost periodic point that was not periodic, so $\overline{\bigcup_{j=1}^\infty \sigma^j a}$ is minimal \cite[p. 457]{ML}.  Since a minimal shift cannot possess any periodic points, this is a contradiction.  Therefore, $S$ cannot contain any other periodic points.
	\end{proof}
	
	\begin{proposition}\label{prop:S-mix}
		$\sigma$ is topologically mixing on $S$.
	\end{proposition}
	\begin{proof}
		To prove $S$ is topologically mixing, it is enough to show it using the language by Definition \ref{def:mixing-words}.  Suppose  $u,v \in \LLL$ such that $u = w_0 0^{a_{i}}\dots w_k0^{a_{i+k}}$ and $v= w'_0 0^{a_{j}}\dots w'_{\ell} 0^{a_{j+\ell}}$.  The case when $i+k < j$ is trivial.  
		
		Suppose $i + k > j$, and note that by assumption \eqref{item:repeat} on the sequence $a$, we know that there exists an $m$ such that $i + k < j + m$ and the word $a_{[j,j+\ell]} = a_{[j+m,j+\ell+m]}$.  Let
		\[w_n^* = 1^n 0^{a_{i+k+1}} 11 0^{a_{i+k+2}} 11 \dots 11 0^{a_{j+m-2}} 11 0^{a_{j+m-1}}\]
		and observe the dependence on $n$.  Thus, we have that for all $n\geq 2$
		\[uw_n^*v = u 1^n 0^{a_{i+k+1}} 11 0^{a_{i+k+2}} 11 \dots 0^{a_{j+m-2}} 11 0^{a_{j+m-1}} w'_0 0^{a_{j+m}} \dots 0^{a_{j+m+\ell}} \in \LLL.\]
		Since $m$ is fixed and only depends on $u$ and $v$, it follows that $S$ is topologically mixing.
	\end{proof}
	
	From here we may use an inverse limit to construct a two-sided shift $X$ in the usual way.  If $X$ contained any periodic orbit other than those of $\overline{1}$ and $\overline{01}$, then $S$ would also contain the one-sided versions of them.  Thus, $X$ only contains these two periodic orbits.  Additionally, the same argument in Proposition \ref{prop:S-mix} shows that $(X,\sigma)$ is topologically mixing.  This shows that there exists a topologically mixing homeomorphism with exactly two periodic orbits.


	\bibliography{Topologically-mixing-suspension-flows}
	\bibliographystyle{amsalpha}	
\end{document}